\documentclass[acmsmall]{acmart}

\usepackage{calc}

\usepackage{subcaption}
\newcommand\subcap[1]{\phantomcaption%
       \caption*{\figurename~\thefigure\thesubfigure: #1}}
\usepackage{stmaryrd,algorithm}
\usepackage{algpseudocode,algorithmicx}

\usepackage{multirow,tabularx}
\usepackage{tikz}
\usetikzlibrary{arrows}
\usetikzlibrary{calc}
\usetikzlibrary{decorations.markings}
\usepackage{ifthen}
\usepackage{xcolor}
\usepackage{doi}
\colorlet{color1}{orange}
\colorlet{color2}{cyan}
\colorlet{color3}{green}
\colorlet{color4}{purple}
\colorlet{highcolor}{magenta}

\newcommand{\highlight}[1]{{\color{highcolor}\textbf{#1}}}

\tikzset{line cap=round,line join=round,>=angle 60,
every picture/.style={line width=1pt}}

\tikzset{->-/.style={postaction={decorate,decoration={
        markings,
        mark=at position .5 with {\arrow[#1]{>}}
      }}}}

\newenvironment{3dmove}[2]{\begin{scope}[shift={(9/4.5*#1,4/4*#1, -#2*12.8/4)}]}{\end{scope}}
\newenvironment{2dmove}[2]{\begin{scope}[shift={(#1*2.9,-#2*2.9)}]}{\end{scope}}

\newenvironment{2dtikzpicture}{\begin{tikzpicture}[x=1cm,y=1cm]}{\end{tikzpicture}}
\newenvironment{3dtikzpicture}{\begin{tikzpicture}[x={(1.126cm,-0.2667cm)},y={(0.6667cm,0.6cm)},z={(0cm,1cm)}]}{\end{tikzpicture}}

\newcommand{\drawaxes}[1][3d]{
\ifthenelse{\equal{#1}{3d}}{
\begin{3dmove}{-0.5}{0}
\draw[->,gray] (0,0,0) -- (.4,0,0) node[anchor=west,inner sep=0] {\footnotesize$x$};
\draw[->,gray] (0,0,0) -- (0,.4,0) node[anchor=south west,inner sep=0] {\footnotesize$y$};
\draw[->,gray] (0,0,0) -- (0,0,.4) node[anchor=south,inner sep=0] {\footnotesize$z$};
\end{3dmove}
}{}
\ifthenelse{\equal{#1}{2din3d}}{
\begin{3dmove}{-0.5}{0}
\draw[->,gray] (0,0,0) -- (.4,0,0) node[anchor=west,inner sep=0] {\footnotesize$x$};
\draw[->,gray] (0,0,0) -- (0,.4,0) node[anchor=south west,inner sep=0] {\footnotesize$y$};
\end{3dmove}
}{}
\ifthenelse{\equal{#1}{2d}}{
\begin{2dmove}{-0.5}{0}
\draw[->,gray] (0,0,0) -- (.4,0,0) node[anchor=west,inner sep=0] {\footnotesize$x$};
\draw[->,gray] (0,0,0) -- (0,.4,0) node[anchor=south,inner sep=0] {\footnotesize$y$};
\end{2dmove}
}{}
}

\usepackage{newfloat}

\usepackage{bm}
\usepackage{cleveref}
\crefname{equation}{}{}
\crefname{figure}{figure}{figures}
\Crefname{figure}{Figure}{Figures}

\newcommand{\DOF}[3][1.5mm]{\raisebox{-0.7mm}{\begin{tikzpicture}
\clip (0,0) circle (#1+1pt);
\draw[#2, fill=white] (0,0) circle (#1) node[black] {\small#3};
\end{tikzpicture}}}

\newcommand{\DOFedge}[2][1.5mm]{\DOF[#1]{color2}{#2}}
\newcommand{\DOFface}[2][1.5mm]{\DOF[#1]{color3}{#2}}

\newcommand{\identmat}{\mathrm{I}}
\newcommand{\refel}{R}
\newcommand{\polyspace}{\mathcal{V}}
\newcommand{\dualbasis}{\mathcal{L}}
\newcommand{\dual}[1]{#1^*}

\usepackage{pythonhighlight}

\newcommand{\dx}[1][x]{\,\mathrm{d}#1}

\renewcommand{\vec}[1]{\bm{#1}}

\newcommand{\baseperm}[2][]{B_{#2}}
\newcommand{\comboperm}[2][]{M_{#2}}

\setcopyright{rightsretained}
\copyrightyear{2022}
\acmYear{2022}

\acmJournal{TOMS}

\begin{document}

\title{Construction of arbitrary order finite element degree-of-freedom maps
      on polygonal and polyhedral cell meshes}

\author{Matthew W. Scroggs}
\email{mws48@cam.ac.uk}
\orcid{0000-0002-4658-2443}
\affiliation{%
  \institution{Department of Engineering, University of Cambridge}
  \streetaddress{Trumpington Street}
  \city{Cambridge}
  \country{United Kingdom}
  \postcode{CB2 1PZ}
}
\author{J{\o}rgen S. Dokken}
\email{jsd55@cam.ac.uk}
\orcid{0000-0001-6489-8858}
\affiliation{%
  \institution{Department of Engineering, University of Cambridge}
  \streetaddress{Trumpington Street}
  \city{Cambridge}
  \country{United Kingdom}
  \postcode{CB2 1PZ}
}
\author{Chris N. Richardson}
\email{chris@bpi.cam.ac.uk}
\orcid{0000-0003-3137-1392}
\affiliation{%
  \institution{BP Institute, University of Cambridge}
  \streetaddress{Bullard Laboratories, Madingley Road}
  \city{Cambridge}
  \country{United Kingdom}
  \postcode{CB3 0EZ}
}
\author{Garth N. Wells}
\email{gnw20@cam.ac.uk}
\orcid{0000-0001-5291-7951}
\affiliation{%
  \institution{Department of Engineering, University of Cambridge}
  \streetaddress{Trumpington Street}
  \city{Cambridge}
  \country{United Kingdom}
  \postcode{CB2 1PZ}
}

\begin{abstract}
  We develop a method for generating degree-of-freedom maps for
  arbitrary order Ciarlet-type finite element spaces for any cell shape.
  The approach is based on the composition of permutations and
  transformations by cell sub-entity. Current approaches to generating
  degree-of-freedom maps for arbitrary order problems typically rely on
  a consistent orientation of cell entities that permits the definition
  of a common local coordinate system on shared edges and faces.
  However, while orientation of a mesh is straightforward for simplex
  cells and is a local operation, it is not a strictly local operation
  for quadrilateral cells and in the case of hexahedral cells not all
  meshes are orientable. The permutation and transformation approach is
  developed for a range of element types, including arbitrary degree
  Lagrange, serendipity, and divergence- and curl-conforming elements,
  and for a range of cell shapes. The approach is local and can be
  applied to cells of any shape, including general polytopes and meshes
  with mixed cell types. A number of examples are presented and the
  developed approach has been implemented in open-source libraries.
\end{abstract}


\maketitle

\section{Introduction}

The construction of conforming finite element spaces requires
neighbouring cells to agree on the layout of shared degrees-of-freedom
(DOFs), and---in the case of vector-valued elements---these cells need
to agree on the orientation of shared edges and faces. These
requirements can be handled straightforwardly for low-order spaces with
DOFs associated with vertices or a single DOF on each edge or face, but
are non-trivial for higher-order finite element spaces. The finite
element method permits arbitrary order approximations, but to exploit
this feature it must be possible to construct appropriate
degree-of-freedom maps and to agree the necessary vector orientations
that provide the required continuity.

To illustrate the orientation issue, \cref{fig:edge_match} shows shared
DOFs on an edge as viewed from neighbouring cells. If the cells do not
agree on the orientation of the edge, the DOFs may be mismatched and the
finite element space will not be conforming. For discontinuous spaces,
similar issues arise with the location of quadrature points on a cell
edge or face. Vector-valued spaces bring additional complexity as the
DOFs have an associated direction that needs to be agreed upon by cells
sharing the DOF. In some cases it suffices to store a common sign. In
other cases, agreement between cells may require linear combinations. To
illustrate, two of the DOFs on a face of a higher-degree N\'ed\'elec
first kind element \cite{Nedelec} on two differently ordered
neighbouring cells are shown in \cref{fig:face_match}. It can be seen
that the directions from the point of view of one cell are a linear
combination of the directions from the point of view of the other cell.

\begin{figure}
\begin{subfigure}{.55\textwidth}
  \centering
  \newcommand{\gapsize}{0.9}
\begin{2dtikzpicture}
\coordinate (gap) at (\gapsize,0);
\clip (-.5,-.2) rectangle ($(0.7,2.2)+(gap)$);
\coordinate (A) at (0,0);
\coordinate (B) at (0.14,2.1);
\coordinate (C) at (-1.2,1.12);

\coordinate (D) at ($(A)+(gap)$);
\coordinate (E) at ($(B)+(gap)$);
\coordinate (F) at (2.1,1.12);

\draw [dashed,gray] (A) -- (D);
\draw [dashed,gray] (B) -- (E);

\draw (A) -- (B) -- (C) -- cycle;

\draw[fill=color1] ($(A)!0.25!(B)$) circle (1.5mm) node {\small0};
\draw[fill=color1] ($(A)!0.5!(B)$) circle (1.5mm) node {\small1};
\draw[fill=color1] ($(A)!0.75!(B)$) circle (1.5mm) node {\small2};

\draw (D) -- (E) -- (F) -- cycle;

\draw[fill=color2] ($(D)!0.25!(E)$) circle (1.5mm) node {\small0};
\draw[fill=color2] ($(D)!0.5!(E)$) circle (1.5mm) node {\small1};
\draw[fill=color2] ($(D)!0.75!(E)$) circle (1.5mm) node {\small2};

\draw[->,shorten >=1.5mm, shorten <=1.5mm] ($(A)!0.25!(B)$) -- ($(D)!0.25!(E)$);
\draw[->,shorten >=1.5mm, shorten <=1.5mm] ($(A)!0.5!(B)$) -- ($(D)!0.5!(E)$);
\draw[->,shorten >=1.5mm, shorten <=1.5mm] ($(A)!0.75!(B)$) -- ($(D)!0.75!(E)$);

\end{2dtikzpicture}
\hspace{10mm}
\begin{2dtikzpicture}[line cap=round,line join=round]
\coordinate (gap) at (\gapsize,0);
\clip (-.5,-.2) rectangle ($(0.7,2.2)+(gap)$);
\coordinate (A) at (0,0);
\coordinate (B) at (0.14,2.1);
\coordinate (C) at (-1.2,1.12);

\coordinate (D) at ($(A)+(gap)$);
\coordinate (E) at ($(B)+(gap)$);
\coordinate (F) at (2.1,1.12);

\draw [dashed,gray] (A) -- (D);
\draw [dashed,gray] (B) -- (E);

\draw (A) -- (B) -- (C) -- cycle;

\draw[fill=color1] ($(A)!0.25!(B)$) circle (1.5mm) node {\small0};
\draw[fill=color1] ($(A)!0.5!(B)$) circle (1.5mm) node {\small1};
\draw[fill=color1] ($(A)!0.75!(B)$) circle (1.5mm) node {\small2};

\draw (D) -- (E) -- (F) -- cycle;

\draw[fill=color2] ($(D)!0.25!(E)$) circle (1.5mm) node {\small2};
\draw[fill=color2] ($(D)!0.5!(E)$) circle (1.5mm) node {\small1};
\draw[fill=color2] ($(D)!0.75!(E)$) circle (1.5mm) node {\small0};

\draw[->,shorten >=1.5mm, shorten <=1.5mm] ($(A)!0.25!(B)$) -- ($(D)!0.75!(E)$);
\draw[->,shorten >=1.5mm, shorten <=1.5mm] ($(A)!0.5!(B)$) -- ($(D)!0.5!(E)$);
\draw[->,shorten >=1.5mm, shorten <=1.5mm] ($(A)!0.75!(B)$) -- ($(D)!0.25!(E)$);

\end{2dtikzpicture}
  \subcap{In the example on the left, the local orderings of the DOFs on
  the shared facet agree. On the right, the orderings do not agree.
  Here, a gap between the two cells has been added to help illustrate
  the orderings on each side.}
  \label{fig:edge_match}
\end{subfigure}\hfill\begin{subfigure}{.42\textwidth}
  \centering
  \input{img/face_match}
  \subcap{The orientations of two DOFs on the faces of a higher-degree
  N\'ed\'elec first kind element on a tetrahedron from the point of
  view of two neighbouring cells. The directions on one face are rotated
  relative to the directions on the other face.}
  \label{fig:face_match}
\end{subfigure}
\end{figure}

It is helpful to have an approach to the consistent ordering and
orientation of DOFs on mesh entities for triangular, quadrilateral,
tetrahedral and hexahedral cells, mixed cell meshes, and for arbitrary
polytopes. The latter case would support less common methods, such as
the virtual element method \cite{beirao:2013} which considers arbitrary
polygonal and polyhedral cells. Before presenting an approach for
general polytope cells, we consider the approaches used in some widely
used finite element libraries.

\subsection{Current approaches}
\label{current-approaches}

Traditional finite element libraries tie the degree of the approximation
space to the degree of the input mesh, with `nodes' for point evaluation
DOFs explicitly provided (and numbered) in the input mesh. This provides
a unique global index for DOFs and makes a DOF map trivial to construct
for specific element types. However, the approach allows only a limited
set of finite element spaces to be constructed, usually low-degree
Lagrange elements, and requires a new mesh to change the degree of the
finite element space.

\begin{figure}
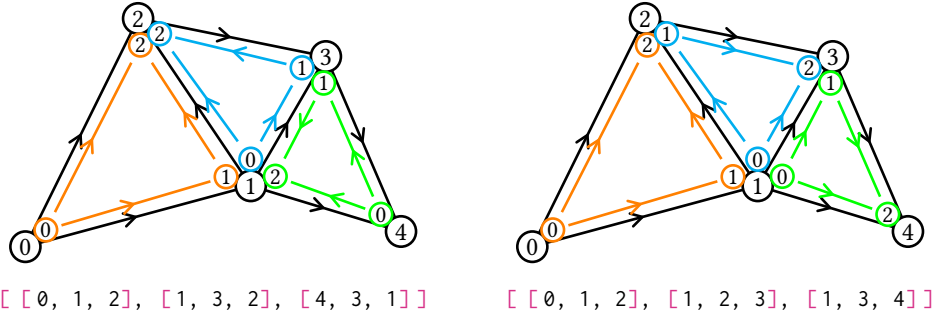

  \newcommand{\localnumber}[5][null]{
\draw[#5, fill=white] ($(#2)+({3.5mm*cos(#3)},{3.5mm*sin(#3)})$) circle (1.5mm) node[black](#1) {\small#4};
}

\begin{2dtikzpicture}
\coordinate (A) at (0,0);
\coordinate (B) at (3,.8);
\coordinate (C) at (1.5,3);
\coordinate (D) at (4,2.5);
\coordinate (E) at (5,0.2);

\draw[->-] (A) -- (B);
\draw[->-] (B) -- (C);
\draw[->-] (A) -- (C);
\draw[->-] (B) -- (D);
\draw[->-] (C) -- (D);
\draw[->-] (B) -- (E);
\draw[->-] (D) -- (E);

\draw[fill=white] (A) circle (2mm) node {0};
\draw[fill=white] (B) circle (2mm) node {1};
\draw[fill=white] (C) circle (2mm) node {2};
\draw[fill=white] (D) circle (2mm) node {3};
\draw[fill=white] (E) circle (2mm) node {4};

\localnumber[t00]{A}{39}{0}{color1}
\localnumber[t01]{B}{160}{1}{color1}
\localnumber[t02]{C}{-85}{2}{color1}

\draw[color1,->-] (t00) -- (t01);
\draw[color1,->-] (t00) -- (t02);
\draw[color1,->-] (t01) -- (t02);

\localnumber[t10]{B}{90}{0}{color2}
\localnumber[t11]{D}{205}{1}{color2}
\localnumber[t12]{C}{-35}{2}{color2}

\draw[color2,->-] (t10) -- (t11);
\draw[color2,->-] (t10) -- (t12);
\draw[color2,->-] (t11) -- (t12);

\localnumber[t20]{E}{140}{0}{color3}
\localnumber[t21]{D}{-95}{1}{color3}
\localnumber[t22]{B}{21}{2}{color3}

\draw[color3,->-] (t20) -- (t21);
\draw[color3,->-] (t20) -- (t22);
\draw[color3,->-] (t21) -- (t22);

\node[anchor=center] at (2.5,-0.7) {\small\pyth{[[0, 1, 2], [1, 3, 2], [4, 3, 1]]}};

\end{2dtikzpicture}
\hspace{5mm}
\begin{2dtikzpicture}
\coordinate (A) at (0,0);
\coordinate (B) at (3,.8);
\coordinate (C) at (1.5,3);
\coordinate (D) at (4,2.5);
\coordinate (E) at (5,0.2);

\draw[->-] (A) -- (B);
\draw[->-] (B) -- (C);
\draw[->-] (A) -- (C);
\draw[->-] (B) -- (D);
\draw[->-] (C) -- (D);
\draw[->-] (B) -- (E);
\draw[->-] (D) -- (E);

\draw[fill=white] (A) circle (2mm) node {0};
\draw[fill=white] (B) circle (2mm) node {1};
\draw[fill=white] (C) circle (2mm) node {2};
\draw[fill=white] (D) circle (2mm) node {3};
\draw[fill=white] (E) circle (2mm) node {4};

\localnumber[t00]{A}{39}{0}{color1}
\localnumber[t01]{B}{160}{1}{color1}
\localnumber[t02]{C}{-85}{2}{color1}

\draw[color1,->-] (t00) -- (t01);
\draw[color1,->-] (t00) -- (t02);
\draw[color1,->-] (t01) -- (t02);

\localnumber[t10]{B}{90}{0}{color2}
\localnumber[t11]{C}{-35}{1}{color2}
\localnumber[t12]{D}{205}{2}{color2}

\draw[color2,->-] (t10) -- (t11);
\draw[color2,->-] (t10) -- (t12);
\draw[color2,->-] (t11) -- (t12);

\localnumber[t20]{B}{21}{0}{color3}
\localnumber[t21]{D}{-95}{1}{color3}
\localnumber[t22]{E}{140}{2}{color3}

\draw[color3,->-] (t20) -- (t21);
\draw[color3,->-] (t20) -- (t22);
\draw[color3,->-] (t21) -- (t22);

\node[anchor=center] at (2.5,-0.7) {\small\pyth{[[0, 1, 2], [1, 2, 3], [1, 3, 4]]}};

\end{2dtikzpicture}
  \caption{An example mesh with three cells. The lists underneath each
  mesh show the global numbers of the vertices of each cell, and the
  coloured numbers show the local ordering of the vertices inside each
  cell. The mesh on the left has been defined so that the vertices of
  each cell are locally ordered in anticlockwise order, but the mesh is
  not ordered, so the local low-to-high edge orientation does not match
  the global low-to-high orientation. The mesh on the right is the
  result of ordering the vertices of each cell of the mesh on the left
  in ascending order: the local order of vertices inside each cell now
  agrees with the global numbering of the vertices, and so the local
  low-to-high edge orientation matches the global low-to-high
  orientation. The vertices, however, are no longer locally ordered in
  an anticlockwise order. If the mesh of triangles is embedded in a 3D
  space---for example if it is a mesh representing the surface of an
  object---this leads to an ambiguity when deciding on the normal vector
  to each face. In the mesh on the left, these normals all naturally
  point out of the page. For the ordered mesh on the right, the face
  normals for the leftmost cell will point out of the page, but for the
  other two cells, it points into the page. For some applications
  the cell orientation still needs to stored despite the cell
  being reordered.}
  \label{fig:ordering}
\end{figure}

Modern finite element libraries usually decouple the input mesh from the
family and degree of the finite element space that is defined on the
mesh. The mesh defines the cell geometry, and different finite element
spaces can then be defined on a given mesh. To construct DOF maps and to
agree on the orientation of any vectors, it is common to require that
the cells in a mesh are locally numbered in a way that guarantees that
entities are consistently ordered. For meshes with simplex cells, it is
straightforward to reorder the local vertex indices of cells to ensure
that the local orientations of the mesh entities on the reference
element agree with their globally agreed orientations; the global
orientation can be inferred from the global vertex indices. The
reordering can be performed locally by cell without reference to
neighbouring cells. An example triangular mesh is shown on the left in
\cref{fig:ordering}. The mesh obtained by sorting the vertices in each
cell in ascending order is shown on the right in \cref{fig:ordering}. In
this reordered mesh, the local and global directions of each edge agree.

The reordering approach is implemented for simplex cells in a number
software libraries, including MFEM \cite{MFEM}, FEniCS
\cite{fenicsbook,alnaes:2015,logg:2010,rognes:2009} and Firedrake
\cite{firedrake}. Once this reordering is done, it is guaranteed that
the orientations of mesh entities are consistent between cells, allowing
a DOF map to be created based on global indices assigned to each mesh
entity, and the directions of any vectors will be consistent between
cells.  NGSolve \cite{NGsolve}---a finite element solver based around
the NetGen mesh generator \cite{NetGen}---avoids having to reorder the
cells by generating meshes in a way that ensures that orientations are
consistent. Higher degree elements are only supported on meshes imported
from an external source when they have a consistent ordering.

\begin{figure}
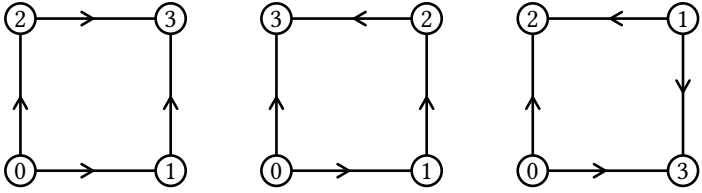

  \begin{2dtikzpicture}
\coordinate (A)  at (0,0);
\coordinate (B)  at (2,0);
\coordinate (C)  at (0,2);
\coordinate (D)  at (2,2);

\draw[->-] (A) -- (C);
\draw[->-] (A) -- (B);
\draw[->-] (C) -- (D);
\draw[->-] (B) -- (D);

\draw[fill=white] (A) circle (2mm) node {0};
\draw[fill=white] (B) circle (2mm) node {1};
\draw[fill=white] (C) circle (2mm) node {2};
\draw[fill=white] (D) circle (2mm) node {3};
\end{2dtikzpicture}
\hspace{8mm}
\begin{2dtikzpicture}
\coordinate (A)  at (0,0);
\coordinate (B)  at (2,0);
\coordinate (C)  at (0,2);
\coordinate (D)  at (2,2);

\draw[->-] (A) -- (C);
\draw[->-] (A) -- (B);
\draw[->-] (D) -- (C);
\draw[->-] (B) -- (D);

\draw[fill=white] (A) circle (2mm) node {0};
\draw[fill=white] (B) circle (2mm) node {1};
\draw[fill=white] (C) circle (2mm) node {3};
\draw[fill=white] (D) circle (2mm) node {2};
\end{2dtikzpicture}
\hspace{8mm}
\begin{2dtikzpicture}
\coordinate (A)  at (0,0);
\coordinate (B)  at (2,0);
\coordinate (C)  at (0,2);
\coordinate (D)  at (2,2);

\draw[->-] (A) -- (C);
\draw[->-] (A) -- (B);
\draw[->-] (D) -- (C);
\draw[->-] (D) -- (B);

\draw[fill=white] (A) circle (2mm) node {0};
\draw[fill=white] (B) circle (2mm) node {3};
\draw[fill=white] (C) circle (2mm) node {2};
\draw[fill=white] (D) circle (2mm) node {1};
\end{2dtikzpicture}
  \caption{These three quadrilateral reference elements are distinct.
  All other local orderings due to the low-to-high convention can be
  obtained by rotating and reflecting these three reference cells.}
  \label{fig:quad_reference}
\end{figure}

The re-ordering approach is a straightforward, local operation for
simplex cells and all meshes of simplex cells can be ordered. For
quadrilateral and hexahedral cells, the reordering approach is more
complicated. \Cref{fig:quad_reference} shows the three distinct
reference cells that would be needed for a naively low-to-high ordered
mesh of quadrilateral cells. For a mesh of hexahedral cells, 501
distinct reference elements would be needed \cite{OEIS-A334304}. deal.II \cite{deal2}
implements a serial re-ordering strategy for quadrilateral and
hexahederal cells to generate consistent edge orderings
\citep{agelek:2017}. After reordering for quadrilateral cells, each cell
is represented by the first reference shown in
\cref{fig:quad_reference}. While re-ordering leads to consistent edge
orientations for all quadrilateral cell meshes of an orientable manifold
(for example, not a M\"obius strip), the operation is not local to each
cell. For the case of hexahedera, not all meshes can be re-ordered to
give a consistent edge direction. If a consistent edge ordering for a
mesh of hexahedera cannot be computed, deal.II will instead store
orientation flags for each subentity of the cell so that these can be
used during assembly \citep{agelek:2017}. A distributed memory parallel
version of the algorithm in deal.II for quadrilaterals is presented
in~\citet{homolya:2016}, and is used in the Firedrake library.

While the re-ordering approach allows high-degree finite element spaces
to be constructed with ease for a range of problems, it is not without
drawbacks. Already mentioned, not all meshes of hexahedral cells can be
reordered. In other cases, explicit orientation markers are still
required: for example, if using a two-dimensional mesh embedded in a
three-dimensional space, or when solving a domain decomposition problem
by coupling the finite element method with the boundary element method
\cite{NedelecFEMBEM,BemppFEMBEM}, the normals to a two-dimensional mesh
are required. These normals are reversed in reflected cells, so whether
or not the cells are reflected must be stored.

When using the $p$-version or $hp$-version of the finite element method
\cite{szabo,ainsworth-coyle,z-thesis} it is common to construct
hierarchical basis functions `on-the-fly' rather than \emph{a priori}
via the Ciarlet definition of a finite element. Basis functions in the
$p$-version are typically associated with cell entities, with `blending
functions' used to ensure that functions associated with a cell entity
have the required trace on other cell entities. As shown in \cite{gatto}
for $H^{1}$ spaces, continuity for higher degree spaces can be easily
guaranteed with this approach by permuting the local coordinates on a
cell entity such that all cells sharing the entity agree on the
coordinate system, with basis functions then computed for each cell.
This approach was extended to some $H({\rm div})$ and $H({\rm curl})$
elements in \cite{fuentes}. An advantage of the approach is simplicity.
Drawbacks are that basis functions cannot be tabulated \emph{a priori}
on a reference cell and must be computed per cell at runtime, and it is
not always clear how more unusual basis functions should be constructed.
The former point limits the application of code generation techniques
for offline basis function evaluation.

\subsection{Proposed approach}

We present an approach for ensuring compatible spaces of arbitrary
degree without creating multiple reference elements, mesh re-ordering or
runtime evaluation of the basis, and which can be applied to meshes with
mixed polytope cells. We consider Ciarlet-type elements, and instead of
reordering each cell, we transform the DOFs on each cell sub-entity to
ensure that the required continuity of finite element fields between
cells is maintained. Given unique global indices for each mesh entity,
the approach is entirely local to each cell. Care is required when
considering integrals over facets, and vector-valued basis functions, as
in each of these cases the orientation of the facet is important and is
no longer guaranteed to agree on each side of the facet.

In contrast to the approaches discussed in \cref{current-approaches},
the approach we propose can be applied to an arbitrary element on an
arbitrary cell type without requiring any specialised implementation for
a specific element or cell. While some parts of our approach---in
particular the permutation of the order of scalar point evaluation
DOFs---are already implemented in other libraries, it does not appear
that the full approach for handling $H({\rm div})$ and $H({\rm curl})$
elements is currently used elsewhere. By representing the
transformations as the product of one (for edges) or two (for faces)
base transformations per sub-entity, we vastly reduce the number of
transformation matrices that need to be stored.

The approach developed in this paper is implemented in the open source
FEniCSx libraries \emph{DOLFINx} and \emph{Basix}. The presented
functionality can be accessed directly from the Basix element library
via C++ and Python interfaces, and its application in practice can be
inspected in the finite element library DOLFINx. The libraries are
available at \url{https://github.com/FEniCS}.

\subsection{Outline}

The remainder of this  paper is arranged as follows. In \cref{sec:dofs},
we look at the definition of a finite element space and the different
types of DOF that we consider. In \cref{sec:ordering}, we describe the
method for transforming the basis functions for each of these DOF types
to achieve a consistent ordering, before presenting some examples of how
this method can be used in \cref{sec:examples}. In
\cref{sec:implementation}, we examine the method's implementation in
FEniCSx. We finish with some concluding remarks in
\cref{sec:conclusion}.

\section{Degrees-of-freedom}
\label{sec:dofs}

Before we look at transforming basis functions and constructing DOF
maps, we define the finite element spaces and DOFs that we will
consider. In general, finite elements can be defined by the
following~\cite{Ciarlet:1978,commonandunusual}.
\begin{definition}[Ciarlet finite element]
A finite element is defined by the triplet
$(\refel,\polyspace,\dualbasis)$, where
\begin{itemize}
  \item $\refel\subset\mathbb{R}^d$ is the reference element, usually a
  polygon or polyhedron;
  \item $\polyspace$ is a finite dimensional polynomial space on $\refel$
  of dimension $n$;
  \item $\dualbasis :=\{l_0,...,l_{n-1}\}$ is a basis of the dual space
  $\dual{\polyspace} : =\{f:\polyspace\to\mathbb{R}\}$.
\end{itemize}
The basis functions $\{\phi_0,\dots,\phi_{n-1}\}$ of the space $\polyspace$ are
defined by
\begin{equation}
  l_i(\phi_j)=\begin{cases}
  1&i=j,\\
  0&i\not=j.
\end{cases}
\label{eq:functional-01}
\end{equation}
\end{definition}

Additionally, a map is defined that maps the basis functions on the
reference element to functions on an arbitrary cell
\cite{rognes:2009,mapping2,mapping3}. The mapping for an element is
chosen so that it preserves important properties of the element's basis
functions.

\begin{figure}
  \input{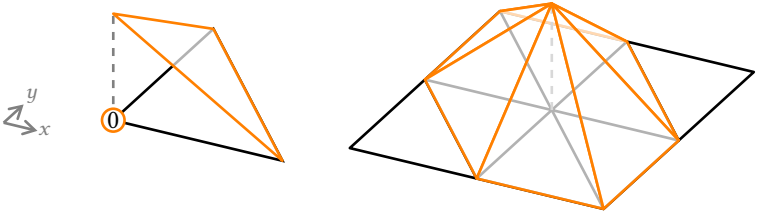}
  \caption{Left: A basis function on the reference of a degree 1
    Lagrange space on a triangle that is associated with vertex 0.
    Right: A basis function of a degree 1 Lagrange space on a triangular
    mesh. Each local DOF associated with the central vertex is
    associated with the same global DOF to ensure continuity.}
  \label{fig:p1_basis_mapped}
\end{figure}

The functionals $l_i \in \dualbasis$ are the (local) DOFs of the finite
element (note that we enumerate functionals and basis functions from~0).
When a finite element is mapped to a mesh, a global DOF number is
assigned to each local DOFs on each cell. To ensure that the mapped
space has the desired continuity properties, each functional $l_i \in
\dualbasis$ is associated with a sub-entity of the reference element. On
neighbouring cells, local DOFs that are associated with a shared
sub-entity are assigned the same global DOF. For example,
\cref{fig:p1_basis_mapped} shows a basis function in a degree 1 Lagrange
space that is associated with vertex 0 and a basis function of a degree
1 Lagrange space on a mesh: the local DOFs on each triangle that are
associated with the central vertex have each been assigned the same
global DOF, causing the function to be continuous.

For clarity, we define the terms we will use to describe meshes and
cells that we will use throughout. By \emph{entities}, we refer to any
of the vertices, edges, faces, or volumes of a cell: these have
dimension 0, 1, 2, and 3 respectively. The \emph{topological dimension}
is the dimension of the cell itself. The \emph{geometric dimension} is
the dimension of the space in which the mesh is embedded. For example, a
mesh of a subset of $\mathbb{R}^2$ with triangular cells has topological
dimension 2 and geometric dimension~2; and a mesh of a 2-dimensional
manifold in $\mathbb{R}^3$ with triangular cells has topological
dimension 2 and geometric dimension~3.
The \emph{codimension} of an entity is given by subtracting the
dimension of the entity from the topological dimension of the cell.
Entities of codimension 1, 2, and 3 are called \emph{facets},
\emph{ridges} and \emph{peaks} (respectively). The usual names given to
entities of cells of topological dimensions 0 to 3 are shown in
\cref{table:entities}.

\begin{table}
\renewcommand{\arraystretch}{1.2}
\begin{tabularx}{\textwidth}{c|cccc|cccc}
\multirow{2}{20mm}{\centering Topological dimension}&
\multicolumn{4}{c|}{Entities by dimension}&
\multicolumn{4}{c}{Entities by codimension}\\
&0&1&2&3&
0&1 (facets)&2 (ridges) &3 (peaks)\\\hline
0 \small(a vertex)
&the cell&-&-&-
&the cell&-&-&-\\
1 \small(an interval)
&vertices&the cell&-&-
&the cell&vertices&-&-\\
2 \small(a polygon)
&vertices&edges&the cell&-
&the cell&edges&vertices&-\\
3 \small(a polyhedron)
&vertices&edges&faces&the cell
&the cell&faces&edges&vertices\\
\end{tabularx}
\vspace{2mm}
\caption{The entities of cells with topological dimensions 0 to~3.}
\label{table:entities}
\end{table}

We now proceed to examine the different types of functional that are
commonly used to define DOFs. The example elements in the following
sections are based on the definitions in the Periodic Table of the
Finite Elements \cite{femtable} and DefElement \cite{defelement}.

\subsection{Point evaluations}
\label{dof:pointeval}

The first---and simplest---functional type is a point evaluation, given
by evaluating the function at a fixed point on the reference element. In
other words, these DOFs are defined, for $v\in\polyspace$, by
\begin{equation}
  l_i(v) := v(\vec{p}_i),
\end{equation}
for some $\vec{p}_i\in\refel$. For example, a degree 1 (linear)
Lagrange finite element on a triangle is defined by taking
$\refel$ to be the triangle with vertices at $(0,0)$, $(1,0)$, and $(0,1)$,
$\polyspace := \operatorname{span}\{1, x, y\}$,
and placing a point evaluation DOF at each of the corners of the
triangle. An example basis function of this space was shown in
\cref{fig:p1_basis_mapped}; note that the basis function is equal to 1
at one vertex and 0 at the other vertices, and so it satisfies
\cref{eq:functional-01}.

\subsection{Point evaluations of face tangents}
\label{dof:pointevaltang}

\begin{figure}
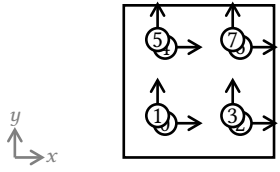

  \begin{2dtikzpicture}
\drawaxes[2d]
\coordinate (A) at (0,0,0);
\coordinate (B) at (2,0,0);
\coordinate (D) at (0,2,0);
\coordinate (C) at (2,2,0);
\draw (A) -- (B) -- (C) -- (D) -- cycle;

\begin{scope}[shift={(0.5,0.5)}]
\draw[black,->] (0.05,-0.05) -- (0.55,-0.05);
\draw[fill=white] (0.05,-0.05) circle (1.5mm) node {\small0};
\draw[black,->] (-0.05,0.05) -- (-0.05,0.55);
\draw[fill=white] (-0.05,0.05) circle (1.5mm) node {\small1};
\end{scope}

\begin{scope}[shift={(1.5,0.5)}]
\draw[black,->] (0.05,-0.05) -- (0.55,-0.05);
\draw[fill=white] (0.05,-0.05) circle (1.5mm) node {\small2};
\draw[black,->] (-0.05,0.05) -- (-0.05,0.55);
\draw[fill=white] (-0.05,0.05) circle (1.5mm) node {\small3};
\end{scope}

\begin{scope}[shift={(0.5,1.5)}]
\draw[black,->] (0.05,-0.05) -- (0.55,-0.05);
\draw[fill=white] (0.05,-0.05) circle (1.5mm) node {\small4};
\draw[black,->] (-0.05,0.05) -- (-0.05,0.55);
\draw[fill=white] (-0.05,0.05) circle (1.5mm) node {\small5};
\end{scope}

\begin{scope}[shift={(1.5,1.5)}]
\draw[black,->] (0.05,-0.05) -- (0.55,-0.05);
\draw[fill=white] (0.05,-0.05) circle (1.5mm) node {\small6};
\draw[black,->] (-0.05,0.05) -- (-0.05,0.55);
\draw[fill=white] (-0.05,0.05) circle (1.5mm) node {\small7};
\end{scope}

\end{2dtikzpicture}
  \caption{A possible arrangment of vector point evaluation DOFs in
  directions tangential to the face.}
  \label{fig:face_tangent_example}
\end{figure}

For vector-valued function spaces, point evaluations must additionally
include a direction in which the evaluation is taken. Vector point
evaluations in directions tangential to a face of a three-dimensional
cell require extra attention. There are two independent vectors tangent
to a face, so it is common to collocate two DOFs at each point.
\Cref{fig:face_tangent_example} shows how eight face tangent DOFs could
be arranged on a quadrilateral face. If two cells meet at a face that
has face tangent DOFs, care must be taken to ensure that the directions
assigned to the DOFs by each cell agree.

\subsection{Integral moments}
\label{dof:intmom}

Numerous finite elements are defined using integral moment DOFs. For
scalar-valued function spaces, these are defined, for $v \in
\polyspace$, by
\begin{equation}
l_i(v) := \int_\Xi v \phi,
\end{equation}
where $\Xi \subseteq \refel$ is some subset of the reference element,
usually the whole cell or one of its faces or edges; and $\phi: \Xi \to
\mathbb{R}$ is a weight function.
For vector-valued function spaces, integral moment DOFs are defined, for
$\vec{v}\in\polyspace$, by
\begin{equation}
l_i(v) := \int_\Xi \vec{v}\cdot\vec\phi.
\end{equation}
In this case, $\vec\phi:\Xi\to\mathbb{R}^d$ is a vector-valued weight
function.

\begin{figure}
\input{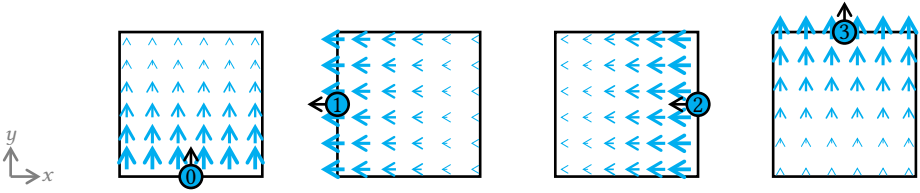}
\caption{The four basis functions of a Raviart--Thomas degree 1 space on
a quadrilateral. The four DOFs are the integral moments over each edge
with the normal to that edge. The normal here have been chosen to be
90{\textdegree} anticlockwise rotations of the low-to-high edge
directions shown in the leftmost reference element in
\cref{fig:quad_reference}.}
\label{fig:quad_div}
\end{figure}

For example, a degree 1 Raviart--Thomas element \cite{RaviartThomas} on
a quadrilateral is defined by taking $\refel$ to be the quadrilateral
with vertices at $(0,0)$, $(1,0)$, $(0,1)$, and $(1,1)$,
\begin{align*}
\polyspace&:=\operatorname{span}\left\{
\begin{pmatrix}1\\0\end{pmatrix},
\begin{pmatrix}x\\0\end{pmatrix},
\begin{pmatrix}0\\1\end{pmatrix},
\begin{pmatrix}0\\y\end{pmatrix}
\right\},
\end{align*}
and defining four DOFs by
\begin{equation}
l_{i}(\vec{v}):=\int_{e_i}\vec{v}\cdot\hat{\vec{n}}_{i}\dx[s],
\end{equation}
where $e_i$ is the $i$th edge of the qudarilateral, and
$\hat{\vec{n}}_{i}$ is the unit vector normal to~$e_i$. The basis
functions of this space are shown in \cref{fig:quad_div}.

As a second example, a degree 2 N\'ed\'elec first kind edge element
\cite{Nedelec} on a tetrahedron is defined by taking $\refel$ to be the
tetrahedron with vertices at $(0,0,0)$, $(1,0,0)$, $(0,1,0)$, and
$(0,0,1)$, and
\begin{align*}
\polyspace&:=\operatorname{span}\left\{
\begin{array}{l}
\begin{pmatrix}1\\0\\0\end{pmatrix},
\begin{pmatrix}x\\0\\0\end{pmatrix},
\begin{pmatrix}y\\0\\0\end{pmatrix},
\begin{pmatrix}0\\1\\0\end{pmatrix},
\begin{pmatrix}0\\x\\0\end{pmatrix},
\begin{pmatrix}0\\y\\0\end{pmatrix},
\begin{pmatrix}0\\0\\1\end{pmatrix},
\begin{pmatrix}0\\0\\x\end{pmatrix},
\begin{pmatrix}0\\0\\y\end{pmatrix},
\\
\begin{pmatrix}z^2\\0\\-xz\end{pmatrix},
\begin{pmatrix}xz\\0\\-x^2\end{pmatrix},
\begin{pmatrix}yz\\0\\-xy\end{pmatrix},
\begin{pmatrix}0\\z^2\\-yz\end{pmatrix},
\begin{pmatrix}0\\yz\\-y^2\end{pmatrix},
\begin{pmatrix}0\\xz\\-xy\end{pmatrix},
\begin{pmatrix}xy\\-x^2\\0\end{pmatrix},
\begin{pmatrix}y^2\\-xy\\0\end{pmatrix}
\end{array}
\right\}.
\end{align*}
Twenty DOFs are used to define this space. The first twelve DOFs are
defined by placing two integral moments on each edge; these are given by
\begin{align*}
l_{2i + j}(\vec{v}):=\int_{e_i}\vec{v}\cdot(\vec{t_i}p_{i,j})\dx[s],
\end{align*}
where $e_i$ is the $i$th edge of the tetrahedron, $p_{i,j}$ is the $j$th
basis function of a degree 1 Lagrange finite element space on $e_i$, and
$\hat{\vec{t}}_{i}$ is the unit vector tangential to $e_i$. We define
the remaining eight DOFs as integral moments on the faces of the
tetrahedron; these are given by
\begin{align*}
l_{6+2j+i}(\vec{v}):=\int_{f_j}\vec{v}\cdot\vec\psi_{j,i}\dx[s],
\end{align*}
where $f_j$ is the $i$th face of the tetrahedron, and $\vec\psi_{j,0}$
and $\vec\psi_{j,1}$ are vector pointing in the directions of two of the
edges of face~$f_j$.

\section{Permuting and transforming cell basis functions}
\label{sec:ordering}

Now that we have examined functionals that can be used to define DOFs,
we look at how the reference basis functions can be adjusted on each
cell to ensure that the required continuity between cells is maintained.
The key idea in this section is to consider basis functions by the mesh
entity with which the corresponding DOF is associated, and then treat these
sets independently.

Let $\{\phi_0,...,\phi_{n-1}\}$ be a set of basis functions, as ordered
on the reference cell. Our general approach is to define a matrix
\(M\in\mathbb{R}^{n\times n}\) for each cell such that the transformed
basis functions for that cell $\{\tilde{\phi}_0, ...,
\tilde{\phi}_{n-1}\}$ are given by
\begin{equation}
\begin{pmatrix}
\tilde{\phi}_0 \\ \vdots \\ \tilde{\phi}_{n-1}
\end{pmatrix}
=
M \begin{pmatrix}
\phi_0 \\ \vdots\\ \phi_{n-1}
\end{pmatrix}.
\end{equation}
The transformation $M$ is computed by comparing the orientation of the
entities of each cell to the orientation of the reference element. The
chosen orientation of the reference element, and the numbering of each
entity, is arbitrary. In all examples presented, we follow the ordering
convention used by Basix \cite{basix}: for simplices, this follows
the UFC ordering convention that orders sub-entities by the numbers of
the vertices that are non-incident to the sub-entity \cite{UFC}. For
other cells, the sub-entities are ordered using the numbers of the
vertices that are incident to the sub-entity.

The method proposed here associates a set of transformations to each
finite element that represent how the element's DOFs are arranged and
how they are affected by reflecting or rotating a sub-entity of the
reference element. Although these transformations are the same for some
elements, they vary greatly for other elements. For example, N\'ed\'elec
second kind elements \cite{Nedelec2} have integral moments with the
basis functions of a Raviart--Thomas space on their facets, meaning the
DOFs on these facets do not have the same arrangement as
the DOFs of a Lagrange element (as shown in the top row of
\cref{fig:base_rotref}). For more exotic elements, the DOF arrangements
can be even more unusual: for example, serendipity spaces
\cite{Arnold_2011} and the divergence- and curl-conforming spaces based
on them \cite{Arnold_2013} are defined for hexahedral cells, but the
DOFs on their faces are integral moments against spaces defined on
triangles. The transformations therefore must be associated with the
element and not the cell type. The operator $M$ is then constructed as a
composition of operations for each element cell entity.

We first consider how spaces of scalar-valued functions are transformed,
before moving on to spaces of vector-valued functions. We consider two
types of vector-valued function space: divergence-conforming spaces and
curl-conforming spaces. These are distinguished by the type of
continuity required between cells: a space is divergence-conforming if
the components normal to the facets of the cell are continuous; a space
is curl-conforming if the components tangent to the facets of the cell
are continuous.

\subsection{Scalar spaces}
\label{sec:scalar-spaces}

For spaces of scalar functions, we must ensure that the DOFs on a mesh
entity appear in the same order when viewed from cells adjacent to the
entity, as shown in \cref{fig:edge_match}. Vertices are trivially
handled as they have no orientation. For two-dimensional cells, the
shared entities of interest are edges; for three-dimensional cell, the
shared entities of interest are edges and faces. When permuting a cell,
the direction of each of its edges may be reversed, and hence order of
the DOFs on each edge may need to be reversed. Faces of a cell can be
rotated and reflected when permuting a cell, and so their DOFs may need
to be rotated or reflected.

\begin{figure}
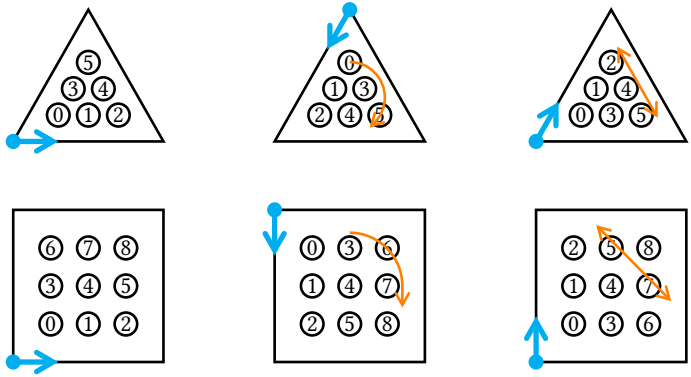

  \begin{2dtikzpicture}
\coordinate (A) at (0,0,0);
\coordinate (B) at (2,0,0);
\coordinate (C) at (1,1.732,0);
\draw (A) -- (B) -- (C) -- cycle;
\draw[color2,->,line width=2pt] (A) -- ($(B)!0.7!(A)$);
\fill[color2] (A) circle (1mm);
\draw[fill=white] ($(A) + ($(A)!0.2!(B)$) + ($(A)!0.2!(C)$)$) circle (1.5mm) node {\small0};
\draw[fill=white] ($(A) + ($(A)!0.4!(B)$) + ($(A)!0.2!(C)$)$) circle (1.5mm) node {\small1};
\draw[fill=white] ($(A) + ($(A)!0.6!(B)$) + ($(A)!0.2!(C)$)$) circle (1.5mm) node {\small2};
\draw[fill=white] ($(A) + ($(A)!0.2!(B)$) + ($(A)!0.4!(C)$)$) circle (1.5mm) node {\small3};
\draw[fill=white] ($(A) + ($(A)!0.4!(B)$) + ($(A)!0.4!(C)$)$) circle (1.5mm) node {\small4};
\draw[fill=white] ($(A) + ($(A)!0.2!(B)$) + ($(A)!0.6!(C)$)$) circle (1.5mm) node {\small5};
\begin{2dmove}{1.2}{0}
\coordinate (A) at (0,0,0);
\coordinate (B) at (2,0,0);
\coordinate (C) at (1,1.732,0);
\draw (A) -- (B) -- (C) -- cycle;
\draw[color2,->,line width=2pt] (C) -- ($(A)!0.7!(C)$);
\fill[color2] (C) circle (1mm);
\draw[fill=white] ($(A) + ($(A)!0.2!(B)$) + ($(A)!0.2!(C)$)$) circle (1.5mm) node {\small2};
\draw[fill=white] ($(A) + ($(A)!0.4!(B)$) + ($(A)!0.2!(C)$)$) circle (1.5mm) node {\small4};
\draw[fill=white] ($(A) + ($(A)!0.6!(B)$) + ($(A)!0.2!(C)$)$) circle (1.5mm) node {\small5};
\draw[fill=white] ($(A) + ($(A)!0.2!(B)$) + ($(A)!0.4!(C)$)$) circle (1.5mm) node {\small1};
\draw[fill=white] ($(A) + ($(A)!0.4!(B)$) + ($(A)!0.4!(C)$)$) circle (1.5mm) node {\small3};
\draw[fill=white] ($(A) + ($(A)!0.2!(B)$) + ($(A)!0.6!(C)$)$) circle (1.5mm) node {\small0};
\draw[color1,->] ($(A)!0.5!(B)!0.3333!(C)+(0,0.4667)$) arc (90:-45:0.4667) -- +(-0.06667,-0.06667);
\end{2dmove}
\begin{2dmove}{2.4}{0}
\coordinate (A) at (0,0,0);
\coordinate (B) at (2,0,0);
\coordinate (C) at (1,1.732,0);
\draw (A) -- (B) -- (C) -- cycle;
\draw[color2,->,line width=2pt] (A) -- ($(C)!0.7!(A)$);
\fill[color2] (A) circle (1mm);
\draw[fill=white] ($(A) + ($(A)!0.2!(B)$) + ($(A)!0.2!(C)$)$) circle (1.5mm) node {\small0};
\draw[fill=white] ($(A) + ($(A)!0.4!(B)$) + ($(A)!0.2!(C)$)$) circle (1.5mm) node {\small3};
\draw[fill=white] ($(A) + ($(A)!0.6!(B)$) + ($(A)!0.2!(C)$)$) circle (1.5mm) node {\small5};
\draw[fill=white] ($(A) + ($(A)!0.2!(B)$) + ($(A)!0.4!(C)$)$) circle (1.5mm) node {\small1};
\draw[fill=white] ($(A) + ($(A)!0.4!(B)$) + ($(A)!0.4!(C)$)$) circle (1.5mm) node {\small4};
\draw[fill=white] ($(A) + ($(A)!0.2!(B)$) + ($(A)!0.6!(C)$)$) circle (1.5mm) node {\small2};
\draw[color1,<->] ($(C)!0.8!(B)!0.1!(A)$) -- ($(B)!0.8!(C)!0.1!(A)$);
\end{2dmove}

\begin{2dmove}{0}{1}
\coordinate (A) at (0,0,0);
\coordinate (B) at (2,0,0);
\coordinate (C) at (2,2,0);
\coordinate (D) at (0,2,0);
\draw (A) -- (B) -- (C) -- (D) -- cycle;
\draw[color2,->,line width=2pt] (A) -- ($(B)!0.7!(A)$);
\fill[color2] (A) circle (1mm);
\draw[fill=white] ($(A) + ($(A)!0.25!(B)$) + ($(A)!0.25!(D)$)$) circle (1.5mm) node {\small0};
\draw[fill=white] ($(A) + ($(A)!0.50!(B)$) + ($(A)!0.25!(D)$)$) circle (1.5mm) node {\small1};
\draw[fill=white] ($(A) + ($(A)!0.75!(B)$) + ($(A)!0.25!(D)$)$) circle (1.5mm) node {\small2};
\draw[fill=white] ($(A) + ($(A)!0.25!(B)$) + ($(A)!0.50!(D)$)$) circle (1.5mm) node {\small3};
\draw[fill=white] ($(A) + ($(A)!0.50!(B)$) + ($(A)!0.50!(D)$)$) circle (1.5mm) node {\small4};
\draw[fill=white] ($(A) + ($(A)!0.75!(B)$) + ($(A)!0.50!(D)$)$) circle (1.5mm) node {\small5};
\draw[fill=white] ($(A) + ($(A)!0.25!(B)$) + ($(A)!0.75!(D)$)$) circle (1.5mm) node {\small6};
\draw[fill=white] ($(A) + ($(A)!0.50!(B)$) + ($(A)!0.75!(D)$)$) circle (1.5mm) node {\small7};
\draw[fill=white] ($(A) + ($(A)!0.75!(B)$) + ($(A)!0.75!(D)$)$) circle (1.5mm) node {\small8};
\end{2dmove}
\begin{2dmove}{1.2}{1}
\coordinate (A) at (0,0,0);
\coordinate (B) at (2,0,0);
\coordinate (C) at (2,2,0);
\coordinate (D) at (0,2,0);
\draw (A) -- (B) -- (C) -- (D) -- cycle;
\draw[color2,->,line width=2pt] (D) -- ($(A)!0.7!(D)$);
\fill[color2] (D) circle (1mm);
\draw[fill=white] ($(A) + ($(A)!0.25!(B)$) + ($(A)!0.25!(D)$)$) circle (1.5mm) node {\small2};
\draw[fill=white] ($(A) + ($(A)!0.50!(B)$) + ($(A)!0.25!(D)$)$) circle (1.5mm) node {\small5};
\draw[fill=white] ($(A) + ($(A)!0.75!(B)$) + ($(A)!0.25!(D)$)$) circle (1.5mm) node {\small8};
\draw[fill=white] ($(A) + ($(A)!0.25!(B)$) + ($(A)!0.50!(D)$)$) circle (1.5mm) node {\small1};
\draw[fill=white] ($(A) + ($(A)!0.50!(B)$) + ($(A)!0.50!(D)$)$) circle (1.5mm) node {\small4};
\draw[fill=white] ($(A) + ($(A)!0.75!(B)$) + ($(A)!0.50!(D)$)$) circle (1.5mm) node {\small7};
\draw[fill=white] ($(A) + ($(A)!0.25!(B)$) + ($(A)!0.75!(D)$)$) circle (1.5mm) node {\small0};
\draw[fill=white] ($(A) + ($(A)!0.50!(B)$) + ($(A)!0.75!(D)$)$) circle (1.5mm) node {\small3};
\draw[fill=white] ($(A) + ($(A)!0.75!(B)$) + ($(A)!0.75!(D)$)$) circle (1.5mm) node {\small6};
\draw[color1,->] ($(A)!0.5!(C)+(0,0.7)$) arc (90:0:0.7) -- +(0,-0.2828);
\end{2dmove}
\begin{2dmove}{2.4}{1}
\coordinate (A) at (0,0,0);
\coordinate (B) at (2,0,0);
\coordinate (C) at (2,2,0);
\coordinate (D) at (0,2,0);
\draw (A) -- (B) -- (C) -- (D) -- cycle;
\draw[color2,->,line width=2pt] (A) -- ($(D)!0.7!(A)$);
\fill[color2] (A) circle (1mm);
\draw[fill=white] ($(A) + ($(A)!0.25!(B)$) + ($(A)!0.25!(D)$)$) circle (1.5mm) node {\small0};
\draw[fill=white] ($(A) + ($(A)!0.50!(B)$) + ($(A)!0.25!(D)$)$) circle (1.5mm) node {\small3};
\draw[fill=white] ($(A) + ($(A)!0.75!(B)$) + ($(A)!0.25!(D)$)$) circle (1.5mm) node {\small6};
\draw[fill=white] ($(A) + ($(A)!0.25!(B)$) + ($(A)!0.50!(D)$)$) circle (1.5mm) node {\small1};
\draw[fill=white] ($(A) + ($(A)!0.50!(B)$) + ($(A)!0.50!(D)$)$) circle (1.5mm) node {\small4};
\draw[fill=white] ($(A) + ($(A)!0.75!(B)$) + ($(A)!0.50!(D)$)$) circle (1.5mm) node {\small7};
\draw[fill=white] ($(A) + ($(A)!0.25!(B)$) + ($(A)!0.75!(D)$)$) circle (1.5mm) node {\small2};
\draw[fill=white] ($(A) + ($(A)!0.50!(B)$) + ($(A)!0.75!(D)$)$) circle (1.5mm) node {\small5};
\draw[fill=white] ($(A) + ($(A)!0.75!(B)$) + ($(A)!0.75!(D)$)$) circle (1.5mm) node {\small8};
\draw[color1,<->] ($(A)!0.5!(C)+(-0.2,0.8)$) -- ($(A)!0.5!(C)+(0.8,-0.2)$);
\end{2dmove}
\end{2dtikzpicture}
  \caption{The effect of rotating and reflecting a face of a tetrahedron
  (top) and hexahedron (bottom) on a possible layout of DOFs on that
  face. The blue arrows have been included to indicate the orientation
  of the triangle or quadrilateral following each transformation.}
  \label{fig:base_rotref}
\end{figure}

DOFs are only associated with a single mesh entity, so permuting one
entity will have no effect on the DOFs of other entities. Each entity
can therefore be considered independently of other entities. We
therefore split the transformation of the DOFs on a cell into \emph{base
transformations}. For each edge of a two- or three-dimensional reference
cell, we store one base transformation that describes the effect of
reversing that edge. For each face of a three-dimensional reference
cell, we store two base transformations: one that describes the effect
of rotating the face; and one that describes the effect of reflecting
the face. The rotations and reflection of a face of a tetrahedron and
hexahedron that are represented by base transformations are shown in
\cref{fig:base_rotref}. For each cell in the mesh, it can then be
calculated how many times each base transformation needs to be applied
to the cell, by comparing the order of the global vertex numbers on each
entity to the ordering of the reference element. We use a low-to-high
orientation of the entities to do this, but other conventions could be
used. The transformation $M$ for the cell can then be calculated by
multiplying each base transformation the necessary number of times.
Matrices representing rotations and reflections of the same face of a
cell do not commute, so it must be decided beforehand in which order
these transformations will appear in the product of base transformations.

For many scalar-valued spaces, the base transformations will be
permutation matrices as permuting a mesh entity leads to swapping the
positions of the DOFs. This is not true in all cases, however: for
higher degree serendipity spaces on a hexahedron, for example, the base
transformations for DOFs on the faces will be linear transformations and
not be permutations. The implementation of higher degree serendipity spaces
will be discussed in \cref{serendipity_example}.

\subsection{Divergence-conforming spaces}
\label{sec:div-spaces}

For vector-valued functions, the direction of the vector used to define
the DOF must be taken into account alongside the ordering of DOFs on an
entity. If two adjacent cells do not agree on the direction of a DOF,
the basis functions on one of the cells
must be adjusted to correct for this.

For divergence-conforming spaces, neighbouring cells must agree on the
direction of the normal vector on facets. If two cells disagree on the
direction of a DOF, then we can multiply the basis functions associated
with the facet on one of the cells by -1 to correct the disagreement. In
order to decide when to multiply values by -1, we need to know whether
or not each entity of a cell has been reflected: this information has
already been calculated when determining how many times to apply each
base transformation to the cell, so this information can be reused here.

\subsection{Curl-conforming spaces}
\label{sec:curl-spaces}

\begin{figure}
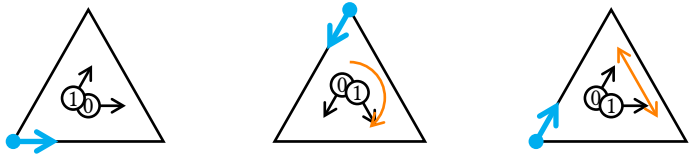

  \begin{2dtikzpicture}
\coordinate (A) at (0,0,0);
\coordinate (B) at (2,0,0);
\coordinate (C) at (1,1.732,0);
\draw (A) -- (B) -- (C) -- cycle;
\draw[color2,->,line width=2pt] (A) -- ($(B)!0.7!(A)$);
\fill[color2] (A) circle (1mm);
\draw[black,->] ($(B)!0.5!(C)!0.4!(A)+(1mm,-0.5mm)$) -- ($(B)!0.5!(C)!0.4!(A)+(1mm,-0.5mm)+(1/2,0)$);
\draw[fill=white] ($(B)!0.5!(C)!0.4!(A)+(1mm,-0.5mm)$) circle (1.5mm) node {\small0};
\draw[black,->] ($(B)!0.5!(C)!0.4!(A)-(1mm,-0.5mm)$) -- ($(B)!0.5!(C)!0.4!(A)-(1mm,-0.5mm)+(0.5/2,0.866/2)$);
\draw[fill=white] ($(B)!0.5!(C)!0.4!(A)-(1mm,-0.5mm)$) circle (1.5mm) node {\small1};
\begin{2dmove}{1.2}{0}
\coordinate (A) at (0,0,0);
\coordinate (B) at (2,0,0);
\coordinate (C) at (1,1.732,0);
\draw (A) -- (B) -- (C) -- cycle;
\draw[color2,->,line width=2pt] (C) -- ($(A)!0.7!(C)$);
\fill[color2] (C) circle (1mm);
\draw[black,->] ($(B)!0.5!(A)!0.4!(C)-(1mm,-0.5mm)$) -- ($(B)!0.5!(A)!0.4!(C)-(1mm,-0.5mm)-(0.5/2,0.866/2)$);
\draw[fill=white] ($(B)!0.5!(A)!0.4!(C)-(1mm,-0.5mm)$) circle (1.5mm) node {\small0};
\draw[black,->] ($(B)!0.5!(A)!0.4!(C)+(1mm,-0.5mm)$) -- ($(B)!0.5!(A)!0.4!(C)+(1mm,-0.5mm)+(0.5/2,-0.866/2)$);
\draw[fill=white] ($(B)!0.5!(A)!0.4!(C)+(1mm,-0.5mm)$) circle (1.5mm) node {\small1};
\draw[color1,->] ($(A)!0.5!(B)!0.3333!(C)+(0,0.4667)$) arc (90:-45:0.4667) -- +(-0.06667,-0.06667);
\end{2dmove}
\begin{2dmove}{2.4}{0}
\coordinate (A) at (0,0,0);
\coordinate (B) at (2,0,0);
\coordinate (C) at (1,1.732,0);
\draw (A) -- (B) -- (C) -- cycle;
\draw[color2,->,line width=2pt] (A) -- ($(C)!0.7!(A)$);
\fill[color2] (A) circle (1mm);
\draw[black,->] ($(B)!0.5!(C)!0.4!(A)-(1mm,-0.5mm)$) -- ($(B)!0.5!(C)!0.4!(A)-(1mm,-0.5mm)+(0.5/2,0.866/2)$);
\draw[fill=white] ($(B)!0.5!(C)!0.4!(A)-(1mm,-0.5mm)$) circle (1.5mm) node {\small0};
\draw[black,->] ($(B)!0.5!(C)!0.4!(A)+(1mm,-0.5mm)$) -- ($(B)!0.5!(C)!0.4!(A)+(1mm,-0.5mm)+(1/2,0)$);
\draw[fill=white] ($(B)!0.5!(C)!0.4!(A)+(1mm,-0.5mm)$) circle (1.5mm) node {\small1};
\draw[color1,<->] ($(C)!0.8!(B)!0.1!(A)$) -- ($(B)!0.8!(C)!0.1!(A)$);
\end{2dmove}
\end{2dtikzpicture}
  \caption{The effect of rotating and reflecting a face of a tetrahedron
  on the direction of face tangents. A rotation leads to $\tilde{l}_0 =
  -l_1$ and $\tilde{l}_1 = l_0-l_1$. A reflection leads to $\tilde{l}_0
  = l_1$ and $\tilde{l}_1 = l_0$. The blue arrows have been included to
  indicate the orientation of the triangle following each
  transformation.}
  \label{fig:face_tangent}
\end{figure}

\begin{figure}
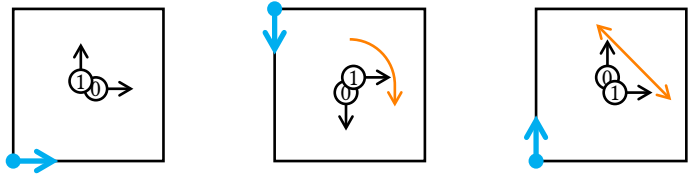

  \begin{2dtikzpicture}
\coordinate (A) at (0,0,0);
\coordinate (B) at (2,0,0);
\coordinate (C) at (2,2,0);
\coordinate (D) at (0,2,0);
\draw (A) -- (B) -- (C) -- (D) -- cycle;
\draw[color2,->,line width=2pt] (A) -- ($(B)!0.7!(A)$);
\fill[color2] (A) circle (1mm);
\draw[black,->] ($(A)!0.5!(C)+(1mm,-0.5mm)$) -- + (1/2,0);
\draw[fill=white] ($(A)!0.5!(C)+(1mm,-0.5mm)$) circle (1.5mm) node {\small0};
\draw[black,->] ($(A)!0.5!(C)-(1mm,-0.5mm)$) -- + (0,1/2);
\draw[fill=white] ($(A)!0.5!(C)-(1mm,-0.5mm)$) circle (1.5mm) node {\small1};
\begin{2dmove}{1.2}{0}
\coordinate (A) at (0,0,0);
\coordinate (B) at (2,0,0);
\coordinate (C) at (2,2,0);
\coordinate (D) at (0,2,0);
\draw (A) -- (B) -- (C) -- (D) -- cycle;
\draw[color2,->,line width=2pt] (D) -- ($(A)!0.7!(D)$);
\fill[color2] (D) circle (1mm);
\draw[black,->] ($(A)!0.5!(C)+(-0.5mm,-1mm)$) -- + (0,-1/2);
\draw[fill=white] ($(A)!0.5!(C)+(-0.5mm,-1mm)$) circle (1.5mm) node {\small0};
\draw[black,->] ($(A)!0.5!(C)-(-0.5mm,-1mm)$) -- + (1/2,0);
\draw[fill=white] ($(A)!0.5!(C)-(-0.5mm,-1mm)$) circle (1.5mm) node {\small1};
\draw[color1,->] ($(A)!0.5!(C)+(0,0.6)$) arc (90:0:0.6) -- +(0,-0.2828);
\end{2dmove}
\begin{2dmove}{2.4}{0}
\coordinate (A) at (0,0,0);
\coordinate (B) at (2,0,0);
\coordinate (C) at (2,2,0);
\coordinate (D) at (0,2,0);
\draw (A) -- (B) -- (C) -- (D) -- cycle;
\draw[color2,->,line width=2pt] (A) -- ($(D)!0.7!(A)$);
\fill[color2] (A) circle (1mm);
\draw[black,->] ($(A)!0.5!(C)+(-0.5mm,1mm)$) -- + (0,1/2);
\draw[fill=white] ($(A)!0.5!(C)+(-0.5mm,1mm)$) circle (1.5mm) node {\small0};
\draw[black,->] ($(A)!0.5!(C)-(-0.5mm,1mm)$) -- + (1/2,0);
\draw[fill=white] ($(A)!0.5!(C)-(-0.5mm,1mm)$) circle (1.5mm) node {\small1};
\draw[color1,<->] ($(A)!0.5!(C)+(-0.2,0.8)$) -- ($(A)!0.5!(C)+(0.8,-0.2)$);
\end{2dmove}
\end{2dtikzpicture}
  \caption{The effect of rotating and reflecting a face of a hexahedron
  on the direction of face tangents. A rotation leads to $\tilde{l}_0 =
  -l_1$ and $\tilde{l}_1 = l_0$. A reflection leads to $\tilde{l}_0 =
  l_1$ and $\tilde{l}_1 = l_0$. The blue arrows have been included to
  indicate the orientation of the quadrilateral following each
  transformation.}
  \label{fig:face_tangent_quad}
\end{figure}

For curl-conforming spaces, the tangential component of the vector will
be continuous between cells. For two-dimensional cells, directions can
be corrected in the same way as for divergence-confirming spaces. For
three-dimensional cells, however, the situation is more complicated. As
we saw in \cref{dof:pointevaltang}, there will commonly be two face
tangent DOFs collocated at points on the face of a three-dimensional
cell, as the space of vectors tangent to the face has dimension~2. For a
given face of a tetrahedral or hexahedral cell, we take the directions
for these DOFs to be parallel to two of the edges of the face (other
directions could be used and would lead to slight variants of the base
transformations we arrive at). The effect of rotating and reflecting the
face on the direction of the DOFs on triangular and quadrilateral faces
are shown in \cref{fig:face_tangent} and \cref{fig:face_tangent_quad}
(respectively).

For face tangent DOFs on faces shared by two cells, linear
transformations must be applied to the basis functions to counteract the
effect of rotating and reflecting a face. As for div-conforming spaces,
the data calculated when working out how many base
transformations to apply to each cell can be reused to determine which
transformations need to be applied.

\subsection{Composing the base transformations for each entity}
\label{sec:composing-transformations}

The transformation matrix $M$ for a given cell can be calculated by
comparing the orientation of each sub-entity to the reference
orientation, then composing (multiplying) the required base
transformations. The transformation matrix can then be applied at cell
level during the finite element assembly process (while the
transformation can be expressed as a matrix multiplication, in practice
it is typically applied without forming $M$).

In cases where the base transformations are permutation matrices, the
permutations can be applied in the construction of the DOF map, avoiding
the application of transformations during assembly. For such
scalar-valued spaces, the permutations we apply are the same as those
proposed in \cite{gatto}, although we apply these permutations to the
DOF numbering of each cell rather than applying them during the creation
of the reference basis functions.

For sub-entities that are edges, we compare the orientation of the edge
to the reference orientation inferred from the global vertex numbers of
the edge's endpoints. If these orientations disagree, then we apply the
base transformation for that edge. Note that in some cases, two cells
adjacent to an edge may agree on the orientation of the edge, but both
disagree with the reference orientation: in this case both cells will
apply their base permutation, as each cell will only compare the edge's
orientation with the reference and not with the other cell.

For sub-entities that are faces, we use two base transformations: one of
these represents moving the vertices one position clockwise (ie a
120{\textdegree} rotation for a equilateral triangle, or a
90{\textdegree} rotation for a square); the other represents a
reflection of the face. The cell must decide how many times this
rotation needs to be applied to the face, and whether or not to apply
the rotation. This can be done by first identifying which of the face's
vertices has the lowest global number: the number of rotations to apply
is calculated by looking at how much the face needs to be rotated to
move this vertex to the position of the locally lowest numbered vertex.
We then look at the two vertices next to the one that we have just
identified: if the lower numbered vertex is not in the same direction as
the corresponding direction for the local vertex numbering, we apply a
reflection to this face. Once the lowest globally numbered vertex and
its lowest numbered neighbour have been identified, we are able to fully
orient the face in a consistent way to neighbouring cells: the global
numbers of the other vertices are irrelevant, and any permutation of the
order of these can be ignored.

\subsection{Discontinuous spaces and facet integrals}
\label{sec:quadrature_permutation}

When using, for example, discontinuous Galerkin schemes, integrals over
shared facets are computed: it is common for these integrals to involve
jump and average operations. This requires the evaluation of basis
functions from both neighbouring cells at a series of quadrature points
on the facet. In order for integrals to be calculated correctly, the two
cells sharing a facet must agree on the orientation of the facet, so
that both cells agree on the location of the quadrature points. On an
unordered mesh, these orientations may not agree, so the quadrature
points must be rotated and reflected to correct for this. Again, the
data calculated when working out how many times to apply the base
transformations can be reused here. The transformation of quadrature
points cannot be amalgamated into the base transformation matrices, so
must be considered separately. The transformation of quadrature points
for facet (and ridge) integrals is essentially the same as the
transformation of DOFs in a space containing point evaluations with the
quadrature points taking the place of the DOF points.

\section{Examples}
\label{sec:examples}

In this section, we discuss how our method can be applied to four
different elements:
a degree 3 Lagrange element on a triangle,
a degree 2 Raviart--Thomas element on a quadrilateral,
a degree 4 Lagrange space on a hexahedron,
and a degree 2 N\'ed\'elec first kind element on a tetrahedron.

\subsection{Computing the base transformations for a cell}

\subsubsection{Example 1: Lagrange on a triangle}

\begin{figure}
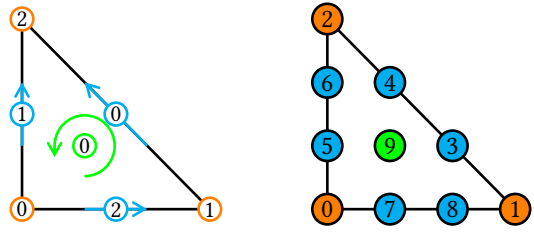

  \begin{2dtikzpicture}
\clip (-.3,-.3) rectangle (2.8,2.8);
\coordinate (A) at (0,0);
\coordinate (B) at (2.5,0);
\coordinate (C) at (0,2.5);

\draw (A) -- (B) -- (C) -- cycle;

\draw[color1, fill=white] (A) circle (1.5mm) node[black] {\small0};
\draw[color1, fill=white] (B) circle (1.5mm) node[black] {\small1};
\draw[color1, fill=white] (C) circle (1.5mm) node[black] {\small2};

\draw[color2, ->] ($(B)!0.3333!(C)$) -- ($(B)!0.6666!(C)$);
\draw[color2, ->] ($(A)!0.3333!(C)$) -- ($(A)!0.6666!(C)$);
\draw[color2, ->] ($(A)!0.3333!(B)$) -- ($(A)!0.6666!(B)$);
\draw[color2, fill=white] ($(B)!0.5!(C)$) circle (1.5mm) node[black] {\small0};
\draw[color2, fill=white] ($(A)!0.5!(C)$) circle (1.5mm) node[black] {\small1};
\draw[color2, fill=white] ($(A)!0.5!(B)$) circle (1.5mm) node[black] {\small2};
\draw[color3, ->] ($(B)!0.5!(C)!0.3333!(A)-(0,4mm)$) arc (-90:180:4mm) -- ($(B)!0.5!(C)!0.3333!(A)-(4mm,1.5mm)$);
\draw[color3, fill=white] ($(B)!0.5!(C)!0.3333!(A)$) circle (1.5mm) node[black] {\small0};
\end{2dtikzpicture}
  \hspace{7mm}
  \begin{2dtikzpicture}
\clip (-.3,-.3) rectangle (2.8, 2.8);
\coordinate (A) at (0,0);
\coordinate (B) at (2.5,0);
\coordinate (C) at (0,2.5);

\draw (A) -- (B) -- (C) -- cycle;

\draw[fill=color1] (A) circle (2mm) node {0};
\draw[fill=color1] (B) circle (2mm) node {1};
\draw[fill=color1] (C) circle (2mm) node {2};
\draw[fill=color2] ($(B)!0.3333!(C)$) circle (2mm) node {3};
\draw[fill=color2] ($(B)!0.6666!(C)$) circle (2mm) node {4};
\draw[fill=color2] ($(A)!0.3333!(C)$) circle (2mm) node {5};
\draw[fill=color2] ($(A)!0.6666!(C)$) circle (2mm) node {6};
\draw[fill=color2] ($(A)!0.3333!(B)$) circle (2mm) node {7};
\draw[fill=color2] ($(A)!0.6666!(B)$) circle (2mm) node {8};
\draw[fill=color3] ($(B)!0.5!(C)!0.3333!(A)$) circle (2mm) node {9};
\end{2dtikzpicture}
  \caption{Left: The numbering and orientation of the entities of the
  reference triangle. Right: The DOF arrangement for a degree 3
  Lagrange space. Each DOF is a point evaluation.}
  \label{fig:p3_reference01}
\end{figure}

The reference element and DOF layout for a degree 3 continuous Lagrange
space are shown in \cref{fig:p3_reference01}. The base transformations
for this element will describe how the order of the DOFs is affected by
reversing an edge. From \cref{fig:p3_reference01}, we can see that
reversing edge 0 will lead to the positions of DOFs 3 and 4 being
swapped; reversing edge 1 will lead to the positions of DOFs 5 and 6
being swapped; and reversing edge 2 will lead to the positions of DOFs 7
and 8 being swapped. These base transformations can be represented by
the following matrices:
\begin{align}
\label{baseperms_eg1}
\baseperm[1]{0}&:=
\begin{pmatrix}
\identmat_3&&&\\
&0&1&\\
&1&0&\\
&&&\identmat_5
\end{pmatrix},
&
\baseperm[1]{1}&:=
\begin{pmatrix}
\identmat_5&&&\\
&0&1&\\
&1&0&\\
&&&\identmat_3
\end{pmatrix},
&
\baseperm[1]{2}&:=
\begin{pmatrix}
\identmat_7&&&\\
&0&1&\\
&1&0&\\
&&&1
\end{pmatrix},
\end{align}
where $\identmat_n$ is an $n$ by $n$ identity matrix and the omitted
blocks of each matrix are~0. For this element, the base transformation
matrices are permutations, so the permutations shown in
\cref{table:p3_base_perms} could be applied directly when creating the
DOF map.

\begin{table}
\begin{tabular}{rrrrrrrrrrll}
0&1&2&\highlight{4}&\highlight{3}&5&6&7&8&9&edge \DOFedge{0} reversed&($\baseperm[1]{0}$)\\[1pt]
0&1&2&3&4&\highlight{6}&\highlight{5}&7&8&9&edge \DOFedge{1} reversed&($\baseperm[1]{1}$)\\[1pt]
0&1&2&3&4&5&6&\highlight{8}&\highlight{7}&9&edge \DOFedge{2} reversed&($\baseperm[1]{2}$)\\[1pt]
\end{tabular}
\vspace{3mm}
\caption{The permutations that the base transformations represent for a
degree 3 Lagrange space on a triangular cell. The DOFs that are on the
entity permuted in each row are highlighted.}
\label{table:p3_base_perms}
\end{table}

\subsubsection{Example 2: Raviart--Thomas on a quadrilateral}

\begin{figure}
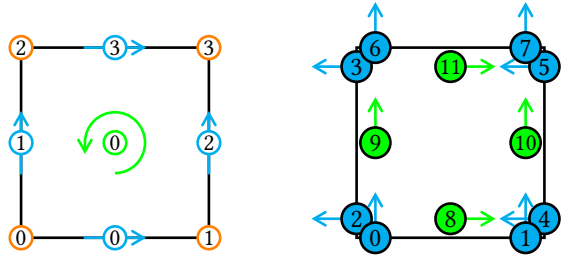

  \begin{2dtikzpicture}
\clip (-.3,-.3) rectangle (2.8,3.2);
\coordinate (A) at (0,0);
\coordinate (B) at (2.5,0);
\coordinate (C) at (0,2.5);
\coordinate (D) at (2.5,2.5);

\draw (A) -- (B) -- (D) -- (C) -- cycle;

\draw[color1, fill=white] (A) circle (1.5mm) node[black] {\small0};
\draw[color1, fill=white] (B) circle (1.5mm) node[black] {\small1};
\draw[color1, fill=white] (C) circle (1.5mm) node[black] {\small2};
\draw[color1, fill=white] (D) circle (1.5mm) node[black] {\small3};

\draw[color2, ->] ($(A)!0.3333!(C)$) -- ($(A)!0.6666!(C)$);
\draw[color2, ->] ($(A)!0.3333!(B)$) -- ($(A)!0.6666!(B)$);
\draw[color2, ->] ($(C)!0.3333!(D)$) -- ($(C)!0.6666!(D)$);
\draw[color2, ->] ($(B)!0.3333!(D)$) -- ($(B)!0.6666!(D)$);
\draw[color2, fill=white] ($(A)!0.5!(B)$) circle (1.5mm) node[black] {\small0};
\draw[color2, fill=white] ($(A)!0.5!(C)$) circle (1.5mm) node[black] {\small1};
\draw[color2, fill=white] ($(B)!0.5!(D)$) circle (1.5mm) node[black] {\small2};
\draw[color2, fill=white] ($(C)!0.5!(D)$) circle (1.5mm) node[black] {\small3};
\draw[color3, ->] ($(A)!0.5!(D)-(0,4mm)$) arc (-90:180:4mm) -- +(0,-1.5mm);
\draw[color3, fill=white] ($(A)!0.5!(D)$) circle (1.5mm) node[black] {\small0};
\end{2dtikzpicture}
  \hspace{7mm}
  \begin{2dtikzpicture}
\clip (-.7,-.3) rectangle (2.8,3.2);
\coordinate (A) at (0,0);
\coordinate (B) at (2.5,0);
\coordinate (C) at (0,2.5);
\coordinate (D) at (2.5,2.5);

\draw (A) -- (B) -- (D) -- (C) -- cycle;

\draw[color2,->] ($(A)!0.1!(C)$) -- +(-0.6,0);
\draw[fill=color2] ($(A)!0.1!(C)$) circle (2mm) node {2};
\draw[color2,->] ($(A)!0.9!(C)$) -- +(-0.6,0);
\draw[fill=color2] ($(A)!0.9!(C)$) circle (2mm) node {3};
\draw[color2,->] ($(B)!0.1!(D)$) -- +(-0.6,0);
\draw[fill=color2] ($(B)!0.1!(D)$) circle (2mm) node {4};
\draw[color2,->] ($(B)!0.9!(D)$) -- +(-0.6,0);
\draw[fill=color2] ($(B)!0.9!(D)$) circle (2mm) node {5};
\draw[color3,->] ($($(A)!0.1!(C)$)!0.5!($(B)!0.1!(D)$)$) -- +(0.6,0);
\draw[fill=color3] ($($(A)!0.1!(C)$)!0.5!($(B)!0.1!(D)$)$) circle (2mm) node {8};
\draw[color3,->] ($($(A)!0.9!(C)$)!0.5!($(B)!0.9!(D)$)$) -- +(0.6,0);
\draw[fill=color3] ($($(A)!0.9!(C)$)!0.5!($(B)!0.9!(D)$)$) circle (2mm) node {\small11};

\draw[color2,->] ($(A)!0.1!(B)$) -- +(0,0.6);
\draw[fill=color2] ($(A)!0.1!(B)$) circle (2mm) node {0};
\draw[color2,->] ($(A)!0.9!(B)$) -- +(0,0.6);
\draw[fill=color2] ($(A)!0.9!(B)$) circle (2mm) node {1};
\draw[color2,->] ($(C)!0.1!(D)$) -- +(0,0.6);
\draw[fill=color2] ($(C)!0.1!(D)$) circle (2mm) node {6};
\draw[color2,->] ($(C)!0.9!(D)$) -- +(0,0.6);
\draw[fill=color2] ($(C)!0.9!(D)$) circle (2mm) node {7};
\draw[color3,->] ($($(A)!0.1!(B)$)!0.5!($(C)!0.1!(D)$)$) -- +(0,0.6);
\draw[fill=color3] ($($(A)!0.1!(B)$)!0.5!($(C)!0.1!(D)$)$) circle (2mm) node {9};
\draw[color3,->] ($($(A)!0.9!(B)$)!0.5!($(C)!0.9!(D)$)$) -- +(0,0.6);
\draw[fill=color3] ($($(A)!0.9!(B)$)!0.5!($(C)!0.9!(D)$)$) circle (2mm) node {\small10};

\end{2dtikzpicture}
  \caption{Left: The numbering and orientation of the entities of the
  reference quadrilateral. Right: The DOF arrangement for a degree 2
  Raviart--Thomas space. DOFs 0 to 7 are integral moments with degree 1
  Lagrange basis functions on the edges. DOFs 8 to 11 are integral
  moments against degree 1 N\'ed\'elec basis functions on the face.}
  \label{fig:rt2_reference01}
\end{figure}

The reference element and DOF layout for a degree 2 Raviart--Thomas
element are shown in \cref{fig:rt2_reference01}. In the same way as in
the first example, the base transformations for this element can be
calculated by considering the effect that reversing an edge would have
on the DOFs on that edge. For this element, however, we must
additionally consider the effect of reversing each edge on the
directions of the DOFs on that edge: when an edge is reversed, these
directions need to be reversed. For example, DOFs 0 and 1 are on edge 0,
so must be reversed if this edge is reversed. Note that DOF 8 is
associated with the interior of the cell, so is unaffected by the
reversal of edge 0. The following base transformation matrices describe
the effect on this element of reversing each edge:
\begin{align}
\baseperm[2]{0}
&:=
\begin{pmatrix}
0&-1&\\
-1&0&\\
&&\identmat_{10}\\
\end{pmatrix},
&
\baseperm[2]{1}
&:=
\begin{pmatrix}
\identmat_2&&&\\
&0&-1&\\
&-1&0&\\
&&&\identmat_8
\end{pmatrix},\notag
\\
\baseperm[2]{2}
&:=
\begin{pmatrix}
\identmat_4&&&\\
&0&-1&\\
&-1&0&\\
&&&\identmat_6
\end{pmatrix},
&
\baseperm[2]{3}
&:=
\begin{pmatrix}
\identmat_6&&&\\
&0&-1&\\
&-1&0&\\
&&&\identmat_4
\end{pmatrix}.
\label{baseperms_eg2}
\end{align}
In each matrix, the DOFs on an edge are swapped (as their positions are
swapped) and multiplied by -1 (as their directions are reversed).

\subsubsection{Example 3: Lagrange on a hexahedron}

\begin{figure}
\input{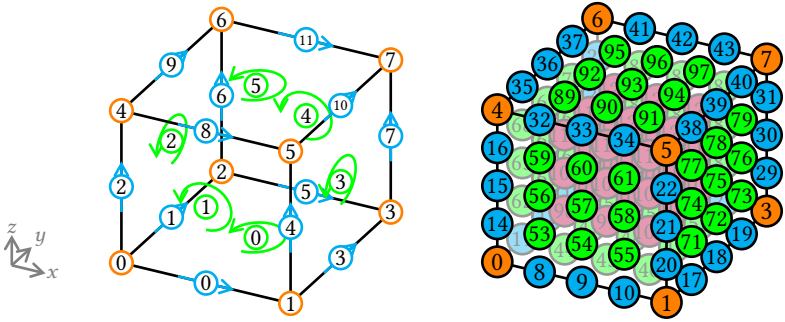}
\hspace{7mm}
\begin{3dtikzpicture}
\newcounter{dofn}
\newcommand{\vdofput}[2]{
\draw[fill=color1] #1 circle (2mm) node {#2};
}
\newcommand{\edofput}[3]{
\setcounter{dofn}{8+3*#3}
\draw[fill=color2] ($#1!0.25!#2$) circle (2mm) node {\ifthenelse{\value{dofn}>9}{\small}{}\arabic{dofn}\addtocounter{dofn}{1}};
\draw[fill=color2] ($#1!0.50!#2$) circle (2mm) node {\ifthenelse{\value{dofn}>9}{\small}{}\arabic{dofn}\addtocounter{dofn}{1}};
\draw[fill=color2] ($#1!0.75!#2$) circle (2mm) node {\ifthenelse{\value{dofn}>9}{\small}{}\arabic{dofn}\addtocounter{dofn}{1}};
}
\newcommand{\fdofput}[6]{
\setcounter{dofn}{8+12*3 + 9*#6}
\draw[fill=color3] ($($#1!0.25!#2$)+($#1!0.25!#3$)-#1$) circle (2mm) node {\ifthenelse{\value{dofn}>9}{\small}{}\arabic{dofn}\addtocounter{dofn}{1}};
\draw[fill=color3] ($($#1!0.50!#2$)+($#1!0.25!#3$)-#1$) circle (2mm) node {\ifthenelse{\value{dofn}>9}{\small}{}\arabic{dofn}\addtocounter{dofn}{1}};
\draw[fill=color3] ($($#1!0.75!#2$)+($#1!0.25!#3$)-#1$) circle (2mm) node {\ifthenelse{\value{dofn}>9}{\small}{}\arabic{dofn}\addtocounter{dofn}{1}};
\draw[fill=color3] ($($#1!0.25!#2$)+($#1!0.50!#3$)-#1$) circle (2mm) node {\ifthenelse{\value{dofn}>9}{\small}{}\arabic{dofn}\addtocounter{dofn}{1}};
\draw[fill=color3] ($($#1!0.50!#2$)+($#1!0.50!#3$)-#1$) circle (2mm) node {\ifthenelse{\value{dofn}>9}{\small}{}\arabic{dofn}\addtocounter{dofn}{1}};
\draw[fill=color3] ($($#1!0.75!#2$)+($#1!0.50!#3$)-#1$) circle (2mm) node {\ifthenelse{\value{dofn}>9}{\small}{}\arabic{dofn}\addtocounter{dofn}{1}};
\draw[fill=color3] ($($#1!0.25!#2$)+($#1!0.75!#3$)-#1$) circle (2mm) node {\ifthenelse{\value{dofn}>9}{\small}{}\arabic{dofn}\addtocounter{dofn}{1}};
\draw[fill=color3] ($($#1!0.50!#2$)+($#1!0.75!#3$)-#1$) circle (2mm) node {\ifthenelse{\value{dofn}>9}{\small}{}\arabic{dofn}\addtocounter{dofn}{1}};
\draw[fill=color3] ($($#1!0.75!#2$)+($#1!0.75!#3$)-#1$) circle (2mm) node {\ifthenelse{\value{dofn}>9}{\small}{}\arabic{dofn}\addtocounter{dofn}{1}};
}

\coordinate (A) at (0,0,0);
\coordinate (B) at (2,0,0);
\coordinate (C) at (0,2,0);
\coordinate (D) at (2,2,0);
\coordinate (E) at (0,0,2);
\coordinate (F) at (2,0,2);
\coordinate (G) at (0,2,2);
\coordinate (H) at (2,2,2);

\draw (G) -- (C) -- (D);
\draw (C) -- (A);

\vdofput{(C)}{2}

\edofput{(A)}{(C)}{1}
\edofput{(C)}{(D)}{5}
\edofput{(C)}{(G)}{6}

\fdofput{(A)}{(B)}{(C)}{(D)}{2}{0}
\fdofput{(C)}{(D)}{(G)}{(H)}{1}{4}

\fdofput{(A)}{(C)}{(E)}{(G)}{0}{2}

\setcounter{dofn}{8+12*3 + 9*6}
\draw[fill=color4] ($($(A)!0.25!(B)$) + ($(A)!0.25!(C)$) + ($(A)!0.25!(E)$)$) circle (2mm) node {\ifthenelse{\value{dofn}>9}{\small}{}\arabic{dofn}\addtocounter{dofn}{1}};
\draw[fill=color4] ($($(A)!0.50!(B)$) + ($(A)!0.25!(C)$) + ($(A)!0.25!(E)$)$) circle (2mm) node {\ifthenelse{\value{dofn}>9}{\small}{}\arabic{dofn}\addtocounter{dofn}{1}};
\draw[fill=color4] ($($(A)!0.75!(B)$) + ($(A)!0.25!(C)$) + ($(A)!0.25!(E)$)$) circle (2mm) node {\ifthenelse{\value{dofn}>9}{\small}{}\arabic{dofn}\addtocounter{dofn}{1}};
\draw[fill=color4] ($($(A)!0.25!(B)$) + ($(A)!0.50!(C)$) + ($(A)!0.25!(E)$)$) circle (2mm) node {\ifthenelse{\value{dofn}>9}{\small}{}\arabic{dofn}\addtocounter{dofn}{1}};
\draw[fill=color4] ($($(A)!0.50!(B)$) + ($(A)!0.50!(C)$) + ($(A)!0.25!(E)$)$) circle (2mm) node {\ifthenelse{\value{dofn}>9}{\small}{}\arabic{dofn}\addtocounter{dofn}{1}};
\draw[fill=color4] ($($(A)!0.75!(B)$) + ($(A)!0.50!(C)$) + ($(A)!0.25!(E)$)$) circle (2mm) node {\ifthenelse{\value{dofn}>9}{\small}{}\arabic{dofn}\addtocounter{dofn}{1}};
\draw[fill=color4] ($($(A)!0.25!(B)$) + ($(A)!0.75!(C)$) + ($(A)!0.25!(E)$)$) circle (2mm) node {\ifthenelse{\value{dofn}>9}{\small}{}\arabic{dofn}\addtocounter{dofn}{1}};
\draw[fill=color4] ($($(A)!0.50!(B)$) + ($(A)!0.75!(C)$) + ($(A)!0.25!(E)$)$) circle (2mm) node {\ifthenelse{\value{dofn}>9}{\small}{}\arabic{dofn}\addtocounter{dofn}{1}};
\draw[fill=color4] ($($(A)!0.75!(B)$) + ($(A)!0.75!(C)$) + ($(A)!0.25!(E)$)$) circle (2mm) node {\ifthenelse{\value{dofn}>9}{\small}{}\arabic{dofn}\addtocounter{dofn}{1}};
\draw[fill=color4] ($($(A)!0.25!(B)$) + ($(A)!0.25!(C)$) + ($(A)!0.50!(E)$)$) circle (2mm) node {\ifthenelse{\value{dofn}>9}{\small}{}\arabic{dofn}\addtocounter{dofn}{1}};
\draw[fill=color4] ($($(A)!0.50!(B)$) + ($(A)!0.25!(C)$) + ($(A)!0.50!(E)$)$) circle (2mm) node {\ifthenelse{\value{dofn}>9}{\small}{}\arabic{dofn}\addtocounter{dofn}{1}};
\draw[fill=color4] ($($(A)!0.75!(B)$) + ($(A)!0.25!(C)$) + ($(A)!0.50!(E)$)$) circle (2mm) node {\ifthenelse{\value{dofn}>9}{\small}{}\arabic{dofn}\addtocounter{dofn}{1}};
\draw[fill=color4] ($($(A)!0.25!(B)$) + ($(A)!0.50!(C)$) + ($(A)!0.50!(E)$)$) circle (2mm) node {\ifthenelse{\value{dofn}>9}{\small}{}\arabic{dofn}\addtocounter{dofn}{1}};
\draw[fill=color4] ($($(A)!0.50!(B)$) + ($(A)!0.50!(C)$) + ($(A)!0.50!(E)$)$) circle (2mm) node {\ifthenelse{\value{dofn}>9}{\small}{}\arabic{dofn}\addtocounter{dofn}{1}};
\draw[fill=color4] ($($(A)!0.75!(B)$) + ($(A)!0.50!(C)$) + ($(A)!0.50!(E)$)$) circle (2mm) node {\ifthenelse{\value{dofn}>9}{\small}{}\arabic{dofn}\addtocounter{dofn}{1}};
\draw[fill=color4] ($($(A)!0.25!(B)$) + ($(A)!0.75!(C)$) + ($(A)!0.50!(E)$)$) circle (2mm) node {\ifthenelse{\value{dofn}>9}{\small}{}\arabic{dofn}\addtocounter{dofn}{1}};
\draw[fill=color4] ($($(A)!0.50!(B)$) + ($(A)!0.75!(C)$) + ($(A)!0.50!(E)$)$) circle (2mm) node {\ifthenelse{\value{dofn}>9}{\small}{}\arabic{dofn}\addtocounter{dofn}{1}};
\draw[fill=color4] ($($(A)!0.75!(B)$) + ($(A)!0.75!(C)$) + ($(A)!0.50!(E)$)$) circle (2mm) node {\ifthenelse{\value{dofn}>9}{\small}{}\arabic{dofn}\addtocounter{dofn}{1}};
\draw[fill=color4] ($($(A)!0.25!(B)$) + ($(A)!0.25!(C)$) + ($(A)!0.75!(E)$)$) circle (2mm) node {\ifthenelse{\value{dofn}>9}{\small}{}\arabic{dofn}\addtocounter{dofn}{1}};
\draw[fill=color4] ($($(A)!0.50!(B)$) + ($(A)!0.25!(C)$) + ($(A)!0.75!(E)$)$) circle (2mm) node {\ifthenelse{\value{dofn}>9}{\small}{}\arabic{dofn}\addtocounter{dofn}{1}};
\draw[fill=color4] ($($(A)!0.75!(B)$) + ($(A)!0.25!(C)$) + ($(A)!0.75!(E)$)$) circle (2mm) node {\ifthenelse{\value{dofn}>9}{\small}{}\arabic{dofn}\addtocounter{dofn}{1}};
\draw[fill=color4] ($($(A)!0.25!(B)$) + ($(A)!0.50!(C)$) + ($(A)!0.75!(E)$)$) circle (2mm) node {\ifthenelse{\value{dofn}>9}{\small}{}\arabic{dofn}\addtocounter{dofn}{1}};
\draw[fill=color4] ($($(A)!0.50!(B)$) + ($(A)!0.50!(C)$) + ($(A)!0.75!(E)$)$) circle (2mm) node {\ifthenelse{\value{dofn}>9}{\small}{}\arabic{dofn}\addtocounter{dofn}{1}};
\draw[fill=color4] ($($(A)!0.75!(B)$) + ($(A)!0.50!(C)$) + ($(A)!0.75!(E)$)$) circle (2mm) node {\ifthenelse{\value{dofn}>9}{\small}{}\arabic{dofn}\addtocounter{dofn}{1}};
\draw[fill=color4] ($($(A)!0.25!(B)$) + ($(A)!0.75!(C)$) + ($(A)!0.75!(E)$)$) circle (2mm) node {\ifthenelse{\value{dofn}>9}{\small}{}\arabic{dofn}\addtocounter{dofn}{1}};
\draw[fill=color4] ($($(A)!0.50!(B)$) + ($(A)!0.75!(C)$) + ($(A)!0.75!(E)$)$) circle (2mm) node {\ifthenelse{\value{dofn}>9}{\small}{}\arabic{dofn}\addtocounter{dofn}{1}};
\draw[fill=color4] ($($(A)!0.75!(B)$) + ($(A)!0.75!(C)$) + ($(A)!0.75!(E)$)$) circle (2mm) node {\ifthenelse{\value{dofn}>9}{\small}{}\arabic{dofn}\addtocounter{dofn}{1}};

\fill[white, opacity=0.6] (A) -- (B) -- (D) -- (H) -- (G) -- (E) -- cycle;

\draw (D) -- (H) -- (G);
\draw (A) -- (B) -- (D);
\draw (B) -- (F) -- (E) -- (A);
\draw (E) -- (G);
\draw (F) -- (H);

\vdofput{(D)}{3}

\vdofput{(A)}{0}
\vdofput{(B)}{1}
\vdofput{(E)}{4}
\vdofput{(F)}{5}
\vdofput{(G)}{6}
\vdofput{(H)}{7}

\edofput{(A)}{(B)}{0}
\edofput{(B)}{(D)}{3}
\edofput{(D)}{(H)}{7}
\edofput{(G)}{(H)}{11}

\edofput{(A)}{(E)}{2}
\edofput{(B)}{(F)}{4}
\edofput{(E)}{(F)}{8}
\edofput{(E)}{(G)}{9}
\edofput{(F)}{(H)}{10}

\fdofput{(A)}{(B)}{(E)}{(F)}{1}{1}
\fdofput{(B)}{(D)}{(F)}{(H)}{0}{3}
\fdofput{(E)}{(F)}{(G)}{(H)}{2}{5}

\end{3dtikzpicture}
\caption{Left: The numbering and orientation of the entities of the
reference hexahedron. Right: The DOF arrangement for a degree 4
Lagrange space. Each DOF is a point evaluation.}
\label{fig:p4_reference01}
\end{figure}

The reference element and DOF layout for a degree 4 Lagrange element
are shown in \cref{fig:p4_reference01}. As in the previous two examples,
the order of the DOFs on the edges of this element will be reversed if
an edge is reversed. As the cell is three dimensional in this example,
we must also include base transformations describing the effect of
rotating and reflecting a face of the cell. Face 1 includes the DOFs 53
to 61. Referring to \cref{fig:base_rotref}, we can see the effect of
rotating this face and construct the base transformations. As in the
first example, each base transformation for this element is a
permutation: these permutations are shown in \cref{table:p4_base_perms}.
We note that the rotation shown in \cref{fig:base_rotref} goes in the
opposite direction to that shown in \cref{fig:p4_reference01}. This is
because we count how many times a face has been rotated to get from a
cell to the reference cell, then use the base transformations to ``undo''
this rotation and map from the reference cell to the physical cell.

\begin{table}
\small
\begin{tabular}{rrrrrrrrrrrll}
\dots&\highlight{10}&\highlight{9}&\highlight{8}&\dots&&&&&&&edge \DOFedge{0} reversed&($\baseperm[3]{0}$)\\[1pt]
\dots&\highlight{13}&\highlight{12}&\highlight{11}&\dots&&&&&&&edge \DOFedge{1} reversed&($\baseperm[3]{1}$)\\[1pt]
\dots&\highlight{16}&\highlight{15}&\highlight{14}&\dots&&&&&&&edge \DOFedge{2} reversed&($\baseperm[3]{2}$)\\[1pt]
\dots&\highlight{19}&\highlight{18}&\highlight{17}&\dots&&&&&&&edge \DOFedge{3} reversed&($\baseperm[3]{3}$)\\[1pt]
\dots&\highlight{22}&\highlight{21}&\highlight{20}&\dots&&&&&&&edge \DOFedge{4} reversed&($\baseperm[3]{4}$)\\[1pt]
\dots&\highlight{25}&\highlight{24}&\highlight{23}&\dots&&&&&&&edge \DOFedge{5} reversed&($\baseperm[3]{5}$)\\[1pt]
\dots&\highlight{28}&\highlight{27}&\highlight{26}&\dots&&&&&&&edge \DOFedge{6} reversed&($\baseperm[3]{6}$)\\[1pt]
\dots&\highlight{31}&\highlight{30}&\highlight{29}&\dots&&&&&&&edge \DOFedge{7} reversed&($\baseperm[3]{7}$)\\[1pt]
\dots&\highlight{34}&\highlight{33}&\highlight{32}&\dots&&&&&&&edge \DOFedge{8} reversed&($\baseperm[3]{8}$)\\[1pt]
\dots&\highlight{37}&\highlight{36}&\highlight{35}&\dots&&&&&&&edge \DOFedge{9} reversed&($\baseperm[3]{9}$)\\[1pt]
\dots&\highlight{40}&\highlight{39}&\highlight{38}&\dots&&&&&&&edge \DOFedge{\tiny10} reversed&($\baseperm[3]{10}$)\\[1pt]
\dots&\highlight{43}&\highlight{42}&\highlight{41}&\dots&&&&&&&edge \DOFedge{\tiny11} reversed&($\baseperm[3]{11}$)\\[1pt]
\dots&\highlight{46}&\highlight{49}&\highlight{52}&\highlight{45}&\highlight{48}&\highlight{51}&\highlight{44}&\highlight{47}&\highlight{50}&\dots&face \DOFface{0} rotated&($\baseperm[3]{12}$)\\[1pt]
\dots&\highlight{44}&\highlight{47}&\highlight{50}&\highlight{45}&\highlight{48}&\highlight{51}&\highlight{46}&\highlight{49}&\highlight{52}&\dots&face \DOFface{0} reflected&($\baseperm[3]{13}$)\\[1pt]
\dots&\highlight{55}&\highlight{58}&\highlight{61}&\highlight{54}&\highlight{57}&\highlight{60}&\highlight{53}&\highlight{56}&\highlight{59}&\dots&face \DOFface{1} rotated&($\baseperm[3]{14}$)\\[1pt]
\dots&\highlight{53}&\highlight{56}&\highlight{59}&\highlight{54}&\highlight{57}&\highlight{60}&\highlight{55}&\highlight{58}&\highlight{61}&\dots&face \DOFface{1} reflected&($\baseperm[3]{15}$)\\[1pt]
\dots&\highlight{64}&\highlight{67}&\highlight{70}&\highlight{63}&\highlight{66}&\highlight{69}&\highlight{62}&\highlight{65}&\highlight{68}&\dots&face \DOFface{2} rotated&($\baseperm[3]{16}$)\\[1pt]
\dots&\highlight{62}&\highlight{65}&\highlight{68}&\highlight{63}&\highlight{66}&\highlight{69}&\highlight{64}&\highlight{67}&\highlight{70}&\dots&face \DOFface{2} reflected&($\baseperm[3]{17}$)\\[1pt]
\dots&\highlight{73}&\highlight{76}&\highlight{79}&\highlight{72}&\highlight{75}&\highlight{78}&\highlight{71}&\highlight{74}&\highlight{77}&\dots&face \DOFface{3} rotated&($\baseperm[3]{18}$)\\[1pt]
\dots&\highlight{71}&\highlight{74}&\highlight{77}&\highlight{72}&\highlight{75}&\highlight{78}&\highlight{73}&\highlight{76}&\highlight{79}&\dots&face \DOFface{3} reflected&($\baseperm[3]{19}$)\\[1pt]
\dots&\highlight{82}&\highlight{85}&\highlight{88}&\highlight{81}&\highlight{84}&\highlight{87}&\highlight{80}&\highlight{83}&\highlight{86}&\dots&face \DOFface{4} rotated&($\baseperm[3]{20}$)\\[1pt]
\dots&\highlight{80}&\highlight{83}&\highlight{86}&\highlight{81}&\highlight{84}&\highlight{87}&\highlight{82}&\highlight{85}&\highlight{88}&\dots&face \DOFface{4} reflected&($\baseperm[3]{21}$)\\[1pt]
\dots&\highlight{91}&\highlight{94}&\highlight{97}&\highlight{90}&\highlight{93}&\highlight{96}&\highlight{89}&\highlight{92}&\highlight{95}&\dots&face \DOFface{5} rotated&($\baseperm[3]{22}$)\\[1pt]
\dots&\highlight{89}&\highlight{92}&\highlight{95}&\highlight{90}&\highlight{93}&\highlight{96}&\highlight{91}&\highlight{94}&\highlight{97}&\dots&face \DOFface{5} reflected&($\baseperm[3]{23}$)
\end{tabular}
\vspace{3mm}
\caption{The permutations described by the base transformations for a
degree 4 Lagrange space on a hexahedral cell. As the space has a
high number of DOFs, only the DOFs on the entity are shown; the other
DOFs will be unaffected by the base transformation.}
\label{table:p4_base_perms}
\end{table}

\subsubsection{Example 4: N\'ed\'elec first kind on a tetrahedron}

\begin{figure}
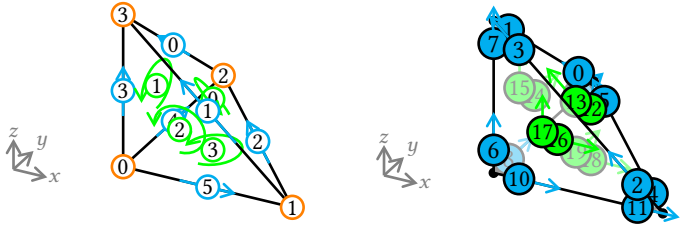

  \begin{3dtikzpicture}
\drawaxes[3d]
\coordinate (A) at (0,0,0);
\coordinate (B) at (2,0,0);
\coordinate (C) at (0,2,0);
\coordinate (D) at (0,0,2);

\draw (A) -- (B) -- (C) -- cycle;

\draw[color2, ->] ($(B)!0.3333!(C)$) -- ($(B)!0.6666!(C)$);
\draw[color2, ->] ($(A)!0.3333!(C)$) -- ($(A)!0.6666!(C)$);
\draw[color2, ->] ($(A)!0.3333!(B)$) -- ($(A)!0.6666!(B)$);
\draw[color2, fill=white] ($(B)!0.5!(C)$) circle (1.5mm) node[black] {\small2};
\draw[color2, fill=white] ($(A)!0.5!(C)$) circle (1.5mm) node[black] {\small4};
\draw[color2, fill=white] ($(A)!0.5!(B)$) circle (1.5mm) node[black] {\small5};
\begin{scope}[shift={($(B)!0.5!(C)!0.3333!(A)$)}]
\draw[color3, ->, domain=-90:180] plot ({0.3*cos(\x)},{0.3*sin(\x)},0) -- (-0.3,-0.05,0);
\draw[color3, fill=white] circle (1.5mm) node[black] {\small3};
\end{scope}

\draw (A) -- (D);
\draw (C) -- (D);

\draw[color1, fill=white] (A) circle (1.5mm) node[black] {\small0};
\draw[color1, fill=white] (C) circle (1.5mm) node[black] {\small2};

\draw[color2, ->] ($(A)!0.3333!(D)$) -- ($(A)!0.6666!(D)$);
\draw[color2, ->] ($(C)!0.3333!(D)$) -- ($(C)!0.6666!(D)$);
\draw[color2, fill=white] ($(A)!0.5!(D)$) circle (1.5mm) node[black] {\small3};
\draw[color2, fill=white] ($(C)!0.5!(D)$) circle (1.5mm) node[black] {\small0};

\begin{scope}[shift={($(D)!0.5!(C)!0.3333!(A)$)}]
\draw[color3, ->, domain=-90:180] plot (0,{0.3*cos(\x)},{0.3*sin(\x)}) -- (0,-0.3,-0.05);
\draw[color3, fill=white] circle (1.5mm) node[black] {\small1};
\end{scope}

\begin{scope}[shift={($(D)!0.5!(B)!0.3333!(C)$)}]
\draw[color3, domain=-90:90] plot ({-0.3/sqrt(2)*cos(\x)-0.15/sqrt(1.5)*sin(\x)},{0.3/sqrt(2)*cos(\x)-0.15/sqrt(1.5)*sin(\x)},{0.3/sqrt(1.5)*sin(\x)});
\draw[color3, fill=white] circle (1.5mm) node[black] {\small0};
\draw[color3, ->, domain=90:180] plot ({-0.3/sqrt(2)*cos(\x)-0.15/sqrt(1.5)*sin(\x)},{0.3/sqrt(2)*cos(\x)-0.15/sqrt(1.5)*sin(\x)},{0.3/sqrt(1.5)*sin(\x)})
        -- ({0.3/sqrt(2)+0.025/sqrt(1.5)},{-0.3/sqrt(2)+0.025/sqrt(1.5)},{-0.05/sqrt(1.5)});
\end{scope}

\draw (B) -- (D);
\draw[color1, fill=white] (B) circle (1.5mm) node[black] {\small1};
\draw[color1, fill=white] (D) circle (1.5mm) node[black] {\small3};

\draw[color2, ->] ($(B)!0.3333!(D)$) -- ($(B)!0.6666!(D)$);
\draw[color2, fill=white] ($(B)!0.5!(D)$) circle (1.5mm) node[black] {\small1};

\begin{scope}[shift={($(D)!0.5!(B)!0.3333!(A)$)}]
\draw[color3, ->, domain=-90:180] plot ({0.3*cos(\x)}, 0,{0.3*sin(\x)}) -- (-0.3,0,-0.05);
\draw[color3, fill=white] circle (1.5mm) node[black] {\small2};
\end{scope}

\end{3dtikzpicture}
  \hspace{7mm}
  \begin{3dtikzpicture}
\drawaxes[3d]
\coordinate (A) at (0,0,0);
\coordinate (B) at (2,0,0);
\coordinate (C) at (0,2,0);
\coordinate (D) at (0,0,2);

\draw (C) -- (A);

\draw[color2,->] ($(A)!0.15!(C)$) -- +(0,0.5,0);
\draw[fill=color2] ($(A)!0.15!(C)$) circle (2mm) node {8};
\draw[color2,->] ($(A)!0.85!(C)$) -- +(0,0.5,0);
\draw[fill=color2] ($(A)!0.85!(C)$) circle (2mm) node {9};

\begin{scope}[shift={($(B)!0.5!(C)!0.3333!(A)$)}]
\draw[color3,->] (1mm,-0.5mm) -- +(0.5,0,0);
\draw[fill=color3] (1mm,-0.5mm) circle (2mm) node {\small18};
\draw[color3,->] (-1mm,0.5mm) -- +(0,0.5,0);
\draw[fill=color3] (-1mm,0.5mm) circle (2mm) node {\small19};
\end{scope}

\begin{scope}[shift={($(D)!0.5!(C)!0.3333!(A)$)}]
\draw[color3,->] (1mm,-0.5mm) -- +(0,0.5,0);
\draw[fill=color3] (1mm,-0.5mm) circle (2mm) node {\small14};
\draw[color3,->] (-1mm,0.5mm) -- +(0,0,0.5);
\draw[fill=color3] (-1mm,0.5mm) circle (2mm) node {\small15};
\end{scope}

\fill[white, opacity=0.6] (B) -- (C) -- (D) -- cycle;
\draw (B) -- (C);
\draw (C) -- (D);

\draw[fill=black] (A) circle (.5mm);
\draw[fill=black] (B) circle (.5mm);
\draw[fill=black] (C) circle (.5mm);
\draw[fill=black] (D) circle (.5mm);

\draw[color2,->] ($(C)!0.15!(D)$) -- +(0,-0.3535,0.3535);
\draw[fill=color2] ($(C)!0.15!(D)$) circle (2mm) node {0};
\draw[color2,->] ($(C)!0.85!(D)$) -- +(0,-0.3535,0.3535);
\draw[fill=color2] ($(C)!0.85!(D)$) circle (2mm) node {1};

\draw[color2,->] ($(B)!0.15!(C)$) -- +(-0.3535,0.3535,0);
\draw[fill=color2] ($(B)!0.15!(C)$) circle (2mm) node {4};
\draw[color2,->] ($(B)!0.85!(C)$) -- +(-0.3535,0.3535,0);
\draw[fill=color2] ($(B)!0.85!(C)$) circle (2mm) node {5};

\begin{scope}[shift={($(D)!0.5!(B)!0.3333!(C)$)}]
\draw[color3,->] (1mm,-0.5mm) -- +(-0.3535,0.3535,0);
\draw[fill=color3] (1mm,-0.5mm) circle (2mm) node {\small12};
\draw[color3,->] (-1mm,0.5mm) -- +(-0.3535,0,0.3535);
\draw[fill=color3] (-1mm,0.5mm) circle (2mm) node {\small13};
\end{scope}

\fill[white, opacity=0.6] (A) -- (B) -- (D) -- cycle;

\draw (A) -- (B) -- (D) -- cycle;

\draw[color2,->] ($(A)!0.15!(D)$) -- +(0,0,0.5);
\draw[fill=color2] ($(A)!0.15!(D)$) circle (2mm) node {6};
\draw[color2,->] ($(A)!0.85!(D)$) -- +(0,0,0.5);
\draw[fill=color2] ($(A)!0.85!(D)$) circle (2mm) node {7};

\draw[color2,->] ($(A)!0.15!(B)$) -- +(0.5,0,0);
\draw[fill=color2] ($(A)!0.15!(B)$) circle (2mm) node {\small10};
\draw[color2,->] ($(A)!0.85!(B)$) -- +(0.5,0,0);
\draw[fill=color2] ($(A)!0.85!(B)$) circle (2mm) node {\small11};

\draw[color2,->] ($(B)!0.15!(D)$) -- +(-0.3535,0,0.3535);
\draw[fill=color2] ($(B)!0.15!(D)$) circle (2mm) node {2};
\draw[color2,->] ($(B)!0.85!(D)$) -- +(-0.3535,0,0.3535);
\draw[fill=color2] ($(B)!0.85!(D)$) circle (2mm) node {3};

\begin{scope}[shift={($(D)!0.5!(B)!0.3333!(A)$)}]
\draw[color3,->] (1mm,-0.5mm) -- +(0.5,0,0);
\draw[fill=color3] (1mm,-0.5mm) circle (2mm) node {\small16};
\draw[color3,->] (-1mm,0.5mm) -- +(0,0,0.5);
\draw[fill=color3] (-1mm,0.5mm) circle (2mm) node {\small17};
\end{scope}

\end{3dtikzpicture}
  \caption{Left: The numbering and orientation of the entities of the
  reference tetrahedron. Right: The DOF arrangement for a degree 2
  N\'ed\'elec first kind curl-conforming space. DOFs 0 to 11 are
  integral moments with degree 1 Lagrange basis functions on the edges.
  DOFs 12 to 19 are integral moments of the tangential components on the
  faces.}
  \label{fig:ned2_reference01}
\end{figure}

The reference element and DOF layout for a degree 2 N\'ed\'elec
first kind element are shown in \cref{fig:ned2_reference01}. Once again,
the effect of reversing each entity can be ascertained by looking at the
layout of DOFs on the reference element. As in the second example, the
directions of DOFs on an edge are reversed when that edge is reversed,
so the base transformations for each edge ($\baseperm[4]{0}$ to
$\baseperm[4]{5}$) will be identity matrices with a 1 replaced with -1
for the DOF on that edge. The DOFs on the faces of this element are not
moved by permuting the face---as the two DOFs on each face are
collocated at the same point---but their directions are affected.
Referring back to \cref{fig:face_tangent}, we can calculate that the
effect of rotating face 0 is described by the transformation
\begin{align*}
\baseperm[4]{6}&:=
\begin{pmatrix}
\identmat_{12}&&&\\
&0&-1&\\
&1&-1&\\
&&&\identmat_6
\end{pmatrix},
\end{align*}
and the effect of reflecting face 0 is described by the transformation
\begin{align*}
\baseperm[4]{7}&:=
\begin{pmatrix}
\identmat_{12}&&&\\
&0&1&\\
&1&0&\\
&&&\identmat_6
\end{pmatrix}.
\end{align*}
These two matrices are the base transformations for face 0, and the
effect of any transformation of face 0 can be calculated by taking
various products of these matrices. We can describe the effect of
rotating and reflecting the other faces in the same way to obtain
$\baseperm[4]{8}$ to~$\baseperm[4]{13}$.

\subsection{Application of the base transformations}

We now look at how the base transformations computed in the previous
section can be applied. \Cref{fig:p3__mesh} shows an example mesh with
four triangular cells that we will compute the transformed basis
functions of a degree 3 Lagrange space on. The reference orientations
of the cell sub-entities are shown in \cref{fig:p3_mesh_numbered_cells}:
the edges marked with a star are those whose local orientation on that
cell do not match the low-to-high orientations shown in
\cref{fig:p3_mesh_def}. These stars tell us which base transformations
we should apply to the reference basis functions for that cell.

\begin{figure}
  \begin{subfigure}{\textwidth}
    \centering
    \input{img/p3_mesh_def}
    \subcap{The mesh of triangles used in the first example, with
    low-to-high edge orientations.}
    \label{fig:p3_mesh_def}
  \end{subfigure}

  \begin{subfigure}{\textwidth}
    \centering
    \input{img/p3_mesh_numbered_cells}
    \subcap{The local numbering and orientations of the entities of the
    four cells in the first and fourth examples. The edges labelled with
    stars are those whose local orientation does not match their global
    orientation shown in \cref{fig:p3_mesh_def}.}
    \label{fig:p3_mesh_numbered_cells}
  \end{subfigure}
  \phantomcaption
  \label{fig:p3__mesh}
\end{figure}

In cell 0, all three edges are reversed, so we should apply all three
base transformations to the reference basis functions. Therefore, for
cell 0 the transformation matrix is given by
\begin{equation}
\comboperm[1]{0}:=
\baseperm[1]{0}\baseperm[1]{1}\baseperm[1]{2}=
\begin{pmatrix}
\identmat_3&&&&&&&\\
&0&1&&&&&\\
&1&0&&&&&\\
&&&0&1&&&\\
&&&1&0&&&\\
&&&&&0&1&\\
&&&&&1&0&\\
&&&&&&&1
\end{pmatrix},
\end{equation}
where $\baseperm[1]{0}$ to $\baseperm[1]{2}$ are the base transformations
defined in \cref{baseperms_eg1}.
In cell 1, no edges are reversed, so $\comboperm[1]{1}$ is the identity
matrix and no change is made to the basis functions. In cell 2, edges
1 and 2 are reversed, so $\comboperm[1]{2}$ is given by
\begin{equation}
\comboperm[1]{2}:=
\baseperm[1]{1}\baseperm[1]{2}=
\begin{pmatrix}
\identmat_5&&&&&\\
&0&1&&&\\
&1&0&&&\\
&&&0&1&\\
&&&1&0&\\
&&&&&1
\end{pmatrix}.
\end{equation}
In cell 4, only edge 0 is reversed, so $\comboperm[1]{3}$ is equal to
$\baseperm[1]{0}$ (the base transformation for edge 0).

As noted above, all the base transformations for this element are
permutations, so they can be applied directly when creating the DOF map.
\Cref{fig:p3_permuted} shows the global DOF numbers assigned to each
cell following the application of these permutations. If the
permutations are not applied, the DOF numbers shown in
\cref{fig:p3_unpermuted} would be assigned instead: note that in this
case, the positions of the DOFs on the edges do not agree from the point
of view of each cell, so the global basis functions would not be
continuous.

\begin{figure}
  \begin{subfigure}{\textwidth}
  \centering
  \newcommand{\nodeDOF}[3]{
\draw[fill=color1] ($(#1)+({3.5mm*cos(#2)},{3.5mm*sin(#2)})$) circle (1.5mm) node {\small#3};
 }

\newcommand{\edgeDOF}[4][]{
\draw[color2,->] ($(#2)!0.3333!(#3)$) -- ($(#2)!0.6666!(#3)$);
\draw[fill=color2] ($(#2)!0.5!(#3)$) circle (1.5mm) node {\small#4};
\ifthenelse{\equal{#1}{}}{}{
\node at ($(#2)!0.5!(#3) + ({3mm*cos(#1)},{3mm*sin(#1)})$) {\color{color2}$\star$};
}
}

\newcommand{\faceDOF}[6][+]{
\ifthenelse{\equal{#1}{+}}{
\draw[color3, ->] ($(#2)!0.5!(#3)!0.3333!(#4)+({4mm*cos(#6)},{4mm*sin(#6)})$)
                  arc (#6:#6+270:4mm)
                  -- +({1.5mm*cos(#6)},{1.5mm*sin(#6)});
}{
\draw[color3, ->] ($(#2)!0.5!(#3)!0.3333!(#4)+({4mm*cos(#6)},{4mm*sin(#6)})$)
                  arc (#6:#6-270:4mm)
                  -- +({1.5mm*cos(#6)},{1.5mm*sin(#6)});
}
\draw[fill=color3] ($(#2)!0.5!(#3)!0.3333!(#4)$) circle (1.5mm) node {\small#5};
}

\begin{2dtikzpicture}
\begin{scope}
\coordinate (p0) at (0,0);
\coordinate (p3) at (1,2);
\coordinate (p4) at (3,2);

\draw (p4) -- (p3) -- (p0) -- cycle;

\coordinate (A) at (p4);
\coordinate (B) at (p3);
\coordinate (C) at (p0);
\draw[fill=color1] (A) circle (2mm) node {\small2};
\draw[fill=color1] (B) circle (2mm) node {\small0};
\draw[fill=color1] (C) circle (2mm) node {\small1};
\draw[fill=color2] ($(B)!0.3333!(C)$) circle (2mm) node {\small4};
\draw[fill=color2] ($(B)!0.6666!(C)$) circle (2mm) node {\small3};
\draw[fill=color2] ($(A)!0.3333!(C)$) circle (2mm) node {\small8};
\draw[fill=color2] ($(A)!0.6666!(C)$) circle (2mm) node {\small7};
\draw[fill=color2] ($(A)!0.3333!(B)$) circle (2mm) node {\small6};
\draw[fill=color2] ($(A)!0.6666!(B)$) circle (2mm) node {\small5};
\draw[fill=color3] ($(B)!0.5!(C)!0.3333!(A)$) circle (2mm) node {\small9};
\node at (0.2,-0.5) {cell 0};
\end{scope}
\begin{scope}[shift={(1.6,0)}]
\coordinate (p0) at (0,0);
\coordinate (p1) at (2,0);
\coordinate (p4) at (3,2);

\draw (p0) -- (p1) -- (p4) -- cycle;

\coordinate (A) at (p0);
\coordinate (B) at (p1);
\coordinate (C) at (p4);
\draw[fill=color1] (A) circle (2mm) node {\small1};
\draw[fill=color1] (B) circle (2mm) node {\small10};
\draw[fill=color1] (C) circle (2mm) node {\small2};
\draw[fill=color2] ($(B)!0.3333!(C)$) circle (2mm) node {\small22};
\draw[fill=color2] ($(B)!0.6666!(C)$) circle (2mm) node {\small19};
\draw[fill=color2] ($(A)!0.3333!(C)$) circle (2mm) node {\small7};
\draw[fill=color2] ($(A)!0.6666!(C)$) circle (2mm) node {\small8};
\draw[fill=color2] ($(A)!0.3333!(B)$) circle (2mm) node {\small20};
\draw[fill=color2] ($(A)!0.6666!(B)$) circle (2mm) node {\small21};
\draw[fill=color3] ($(B)!0.5!(C)!0.3333!(A)$) circle (2mm) node {\small23};
\node at (1,-0.5) {cell 1};
\end{scope}
\begin{scope}[shift={(2.6,0)}]
\coordinate (p1) at (2,0);
\coordinate (p2) at (4,0);
\coordinate (p4) at (3,2);

\draw (p4) -- (p1) -- (p2) -- cycle;

\coordinate (A) at (p4);
\coordinate (B) at (p1);
\coordinate (C) at (p2);
\draw[fill=color1] (A) circle (2mm) node {\small2};
\draw[fill=color1] (B) circle (2mm) node {\small10};
\draw[fill=color1] (C) circle (2mm) node {\small11};
\draw[fill=color2] ($(B)!0.3333!(C)$) circle (2mm) node {\small17};
\draw[fill=color2] ($(B)!0.6666!(C)$) circle (2mm) node {\small18};
\draw[fill=color2] ($(A)!0.3333!(C)$) circle (2mm) node {\small16};
\draw[fill=color2] ($(A)!0.6666!(C)$) circle (2mm) node {\small15};
\draw[fill=color2] ($(A)!0.3333!(B)$) circle (2mm) node {\small19};
\draw[fill=color2] ($(A)!0.6666!(B)$) circle (2mm) node {\small22};
\draw[fill=color3] ($(B)!0.5!(C)!0.3333!(A)$) circle (2mm) node {\small27};
\node at (3,-0.5) {cell 2};
\end{scope}
\begin{scope}[shift={(3.6,0)}]
\coordinate (p2) at (4,0);
\coordinate (p4) at (3,2);
\coordinate (p5) at (5,2);

\draw (p2) -- (p5) -- (p4) -- cycle;

\coordinate (A) at (p2);
\coordinate (B) at (p5);
\coordinate (C) at (p4);
\draw[fill=color1] (A) circle (2mm) node {\small11};
\draw[fill=color1] (B) circle (2mm) node {\small12};
\draw[fill=color1] (C) circle (2mm) node {\small2};
\draw[fill=color2] ($(B)!0.3333!(C)$) circle (2mm) node {\small25};
\draw[fill=color2] ($(B)!0.6666!(C)$) circle (2mm) node {\small24};
\draw[fill=color2] ($(A)!0.3333!(C)$) circle (2mm) node {\small15};
\draw[fill=color2] ($(A)!0.6666!(C)$) circle (2mm) node {\small16};
\draw[fill=color2] ($(A)!0.3333!(B)$) circle (2mm) node {\small13};
\draw[fill=color2] ($(A)!0.6666!(B)$) circle (2mm) node {\small14};
\draw[fill=color3] ($(B)!0.5!(C)!0.3333!(A)$) circle (2mm) node {\small26};

\node at (4,-0.5) {cell 3};
\end{scope}
\end{2dtikzpicture}
  \subcap{The DOF numbering on each cell in the mesh once the
  permutation has been applied to each cell.}
  \label{fig:p3_permuted}
\end{subfigure}

\vspace{2mm}

\begin{subfigure}{\textwidth}
  \centering
  \newcommand{\nodeDOF}[3]{
\draw[fill=color1] ($(#1)+({3.5mm*cos(#2)},{3.5mm*sin(#2)})$) circle (1.5mm) node {\small#3};
 }

\newcommand{\edgeDOF}[4][]{
\draw[color2,->] ($(#2)!0.3333!(#3)$) -- ($(#2)!0.6666!(#3)$);
\draw[fill=color2] ($(#2)!0.5!(#3)$) circle (1.5mm) node {\small#4};
\ifthenelse{\equal{#1}{}}{}{
\node at ($(#2)!0.5!(#3) + ({3mm*cos(#1)},{3mm*sin(#1)})$) {\color{color2}$\star$};
}
}

\newcommand{\faceDOF}[6][+]{
\ifthenelse{\equal{#1}{+}}{
\draw[color3, ->] ($(#2)!0.5!(#3)!0.3333!(#4)+({4mm*cos(#6)},{4mm*sin(#6)})$)
                  arc (#6:#6+270:4mm)
                  -- +({1.5mm*cos(#6)},{1.5mm*sin(#6)});
}{
\draw[color3, ->] ($(#2)!0.5!(#3)!0.3333!(#4)+({4mm*cos(#6)},{4mm*sin(#6)})$)
                  arc (#6:#6-270:4mm)
                  -- +({1.5mm*cos(#6)},{1.5mm*sin(#6)});
}
\draw[fill=color3] ($(#2)!0.5!(#3)!0.3333!(#4)$) circle (1.5mm) node {\small#5};
}

\begin{2dtikzpicture}
\begin{scope}
\coordinate (p0) at (0,0);
\coordinate (p3) at (1,2);
\coordinate (p4) at (3,2);

\draw (p4) -- (p3) -- (p0) -- cycle;

\coordinate (A) at (p4);
\coordinate (B) at (p3);
\coordinate (C) at (p0);
\draw[fill=color1] (A) circle (2mm) node {\small2};
\draw[fill=color1] (B) circle (2mm) node {\small0};
\draw[fill=color1] (C) circle (2mm) node {\small1};
\draw[fill=color2] ($(B)!0.3333!(C)$) circle (2mm) node {\small3};
\draw[fill=color2] ($(B)!0.6666!(C)$) circle (2mm) node {\small4};
\draw[fill=color2] ($(A)!0.3333!(C)$) circle (2mm) node {\small7};
\draw[fill=color2] ($(A)!0.6666!(C)$) circle (2mm) node {\small8};
\draw[fill=color2] ($(A)!0.3333!(B)$) circle (2mm) node {\small5};
\draw[fill=color2] ($(A)!0.6666!(B)$) circle (2mm) node {\small6};
\draw[fill=color3] ($(B)!0.5!(C)!0.3333!(A)$) circle (2mm) node {\small9};
\node at (0.2,-0.5) {cell 0};
\end{scope}
\begin{scope}[shift={(1.6,0)}]
\coordinate (p0) at (0,0);
\coordinate (p1) at (2,0);
\coordinate (p4) at (3,2);

\draw (p0) -- (p1) -- (p4) -- cycle;

\coordinate (A) at (p0);
\coordinate (B) at (p1);
\coordinate (C) at (p4);
\draw[fill=color1] (A) circle (2mm) node {\small1};
\draw[fill=color1] (B) circle (2mm) node {\small10};
\draw[fill=color1] (C) circle (2mm) node {\small2};
\draw[fill=color2] ($(B)!0.3333!(C)$) circle (2mm) node {\small22};
\draw[fill=color2] ($(B)!0.6666!(C)$) circle (2mm) node {\small19};
\draw[fill=color2] ($(A)!0.3333!(C)$) circle (2mm) node {\small7};
\draw[fill=color2] ($(A)!0.6666!(C)$) circle (2mm) node {\small8};
\draw[fill=color2] ($(A)!0.3333!(B)$) circle (2mm) node {\small20};
\draw[fill=color2] ($(A)!0.6666!(B)$) circle (2mm) node {\small21};
\draw[fill=color3] ($(B)!0.5!(C)!0.3333!(A)$) circle (2mm) node {\small23};
\node at (1,-0.5) {cell 1};
\end{scope}
\begin{scope}[shift={(2.6,0)}]
\coordinate (p1) at (2,0);
\coordinate (p2) at (4,0);
\coordinate (p4) at (3,2);

\draw (p4) -- (p1) -- (p2) -- cycle;

\coordinate (A) at (p4);
\coordinate (B) at (p1);
\coordinate (C) at (p2);
\draw[fill=color1] (A) circle (2mm) node {\small2};
\draw[fill=color1] (B) circle (2mm) node {\small10};
\draw[fill=color1] (C) circle (2mm) node {\small11};
\draw[fill=color2] ($(B)!0.3333!(C)$) circle (2mm) node {\small17};
\draw[fill=color2] ($(B)!0.6666!(C)$) circle (2mm) node {\small18};
\draw[fill=color2] ($(A)!0.3333!(C)$) circle (2mm) node {\small15};
\draw[fill=color2] ($(A)!0.6666!(C)$) circle (2mm) node {\small16};
\draw[fill=color2] ($(A)!0.3333!(B)$) circle (2mm) node {\small22};
\draw[fill=color2] ($(A)!0.6666!(B)$) circle (2mm) node {\small19};
\draw[fill=color3] ($(B)!0.5!(C)!0.3333!(A)$) circle (2mm) node {\small27};
\node at (3,-0.5) {cell 2};
\end{scope}
\begin{scope}[shift={(3.6,0)}]
\coordinate (p2) at (4,0);
\coordinate (p4) at (3,2);
\coordinate (p5) at (5,2);

\draw (p2) -- (p5) -- (p4) -- cycle;

\coordinate (A) at (p2);
\coordinate (B) at (p5);
\coordinate (C) at (p4);
\draw[fill=color1] (A) circle (2mm) node {\small11};
\draw[fill=color1] (B) circle (2mm) node {\small12};
\draw[fill=color1] (C) circle (2mm) node {\small2};
\draw[fill=color2] ($(B)!0.3333!(C)$) circle (2mm) node {\small24};
\draw[fill=color2] ($(B)!0.6666!(C)$) circle (2mm) node {\small25};
\draw[fill=color2] ($(A)!0.3333!(C)$) circle (2mm) node {\small15};
\draw[fill=color2] ($(A)!0.6666!(C)$) circle (2mm) node {\small16};
\draw[fill=color2] ($(A)!0.3333!(B)$) circle (2mm) node {\small13};
\draw[fill=color2] ($(A)!0.6666!(B)$) circle (2mm) node {\small14};
\draw[fill=color3] ($(B)!0.5!(C)!0.3333!(A)$) circle (2mm) node {\small26};

\node at (4,-0.5) {cell 3};
\end{scope}
\end{2dtikzpicture}
  \subcap{The DOF numbering on each cell if no permutation is performed.}
  \label{fig:p3_unpermuted}
  \end{subfigure}
\end{figure}

For the degree 2 Raviart--Thomas space on a quadrilateral, the
transformation matrix for each cell can be computed in the same way by
combining the base transformation defined in \cref{baseperms_eg2},
except in this case, the transformations are not permutations, so cannot
be applied directly to the DOF map. For example, for the mesh and
reference orientations shown in \cref{fig:rt2__mesh}, $\comboperm[2]{0}$
would be the identity for cell~0. For cell~1,
$\comboperm[2]{1}:=\baseperm[2]{0}\baseperm[2]{2}$ (the product of the
base permutations for edges 0 and 2).

\begin{figure}
  \begin{subfigure}{0.38\textwidth}
    \centering
    \newcommand{\nodeDOF}[3]{
\draw[fill=color1] ($(#1)+({3.5mm*cos(#2)},{3.5mm*sin(#2)})$) circle (1.5mm) node {\small#3};
 }

\newcommand{\edgeDOF}[4][]{
\draw[color2,->] ($(#2)!0.3333!(#3)$) -- ($(#2)!0.6666!(#3)$);
\draw[fill=color2] ($(#2)!0.5!(#3)$) circle (1.5mm) node {\small#4};
\ifthenelse{\equal{#1}{}}{}{
\node at ($(#2)!0.5!(#3) + ({3mm*cos(#1)},{3mm*sin(#1)})$) {\color{color2}$\star$};
}
}

\newcommand{\faceDOF}[6][+]{
\draw[fill=color3] ($(#2)!0.5!(#3)!0.3333!(#4)$) circle (1.5mm) node {\small#5};
}

\begin{2dtikzpicture}
\coordinate (p0) at (0,0);
\coordinate (p1) at (2,0);
\coordinate (p2) at (4,0);
\coordinate (p3) at (0,2);
\coordinate (p4) at (2,2);
\coordinate (p5) at (4,2);

\draw[->-] (p0) -- (p1);
\draw[->-] (p1) -- (p4);
\draw[->-] (p3) -- (p4);
\draw[->-] (p0) -- (p3);
\draw[->-] (p1) -- (p2);
\draw[->-] (p2) -- (p5);
\draw[->-] (p4) -- (p5);

\draw[fill=white] (p0) circle (2mm) node {0};
\draw[fill=white] (p1) circle (2mm) node {1};
\draw[fill=white] (p2) circle (2mm) node {2};
\draw[fill=white] (p3) circle (2mm) node {3};
\draw[fill=white] (p4) circle (2mm) node {4};
\draw[fill=white] (p5) circle (2mm) node {5};

\node[align=center] at ($($(p0)!0.5!(p4)$)+(0,-0.3)$) {cell 0\\[-3pt]\footnotesize{\pyth{[0, 1, 3, 4]}}};
\node[align=center] at ($($(p4)!0.5!(p2)$)+(0,-0.3)$) {cell 1\\[-3pt]\footnotesize{\pyth{[4, 1, 5, 2]}}};
\end{2dtikzpicture}
    \subcap{The mesh of quadrilaterals used in the second example, with
    low-to-high edge orientations.}
    \label{fig:rt2_mesh_def}
  \end{subfigure}\hfill\begin{subfigure}{0.58\textwidth}
    \centering
    \input{img/rt2_mesh_numbered_cells}
    \subcap{The local numbering and orientations of the entities of the
    three cells in the second example. The edges labelled with stars are
    those whose local orientation does not match their global
    orientation shown in \cref{fig:rt2_mesh_def}.}
    \label{fig:rt2_mesh_numbered_cells}
  \end{subfigure}
  \phantomcaption
  \label{fig:rt2__mesh}
\end{figure}

For three-dimensional cells, the situation is slightly more complex, as
we have to apply base transformations for the faces as well as the
edges. For each face we need to ascertain how many times the face needs
to be rotated and if it needs to be reflected to map from the reference
orientation. A two cell mesh for a degree 4 Lagrange space on a
hexahedron is shown in \cref{fig:p4__mesh}. \Cref{table:p4_reflect_data}
summarises which base transformations should be applied to reference
basis functions for each cell of this mesh. While the transformations
due to reflections can only be applied once (as reflecting twice is the
same as doing nothing), the transformations due to rotations can be
applied more than once: for example, the transformation due to rotating
face~0 must be applied twice to the reference basis functions for
cell~0.

\begin{figure}
  \begin{subfigure}{\textwidth}
  \centering
  \newcommand{\nodeDOF}[3]{
\draw[fill=color1] ($(#1)+({3.5mm*cos(#2)},{3.5mm*sin(#2)})$) circle (1.5mm) node {\small#3};
 }

\newcommand{\edgeDOF}[4][]{
\draw[color2,->] ($(#2)!0.3333!(#3)$) -- ($(#2)!0.6666!(#3)$);
\draw[fill=color2] ($(#2)!0.5!(#3)$) circle (1.5mm) node {\small#4};
\ifthenelse{\equal{#1}{}}{}{
\node at ($(#2)!0.5!(#3) + ({3mm*cos(#1)},{3mm*sin(#1)})$) {\color{color2}$\star$};
}
}

\newcommand{\faceDOF}[6][+]{
\draw[fill=color3] ($(#2)!0.5!(#3)!0.3333!(#4)$) circle (1.5mm) node {\small#5};
}

\begin{3dtikzpicture}
\drawaxes[3d]
\coordinate (p0) at (0,0,0);
\coordinate (p1) at (2.5,0,0);
\coordinate (p2) at (0,2.5,0);
\coordinate (p3) at (2.5,2.5,2.5);
\coordinate (p4) at (0,0,2.5);
\coordinate (p5) at (2.5,2.5,0);
\coordinate (p6) at (2.5,0,2.5);
\coordinate (p7) at (0,2.5,2.5);
\coordinate (p8) at (5,2.5,2.5);
\coordinate (p9) at (5,2.5,0);
\coordinate (p10) at (5,0,2.5);
\coordinate (p11) at (5,0,0);

\draw[->-] (p0) -- (p2);
\draw[->-] (p0) -- (p1);
\draw[->-] (p2) -- (p5);
\draw[->-] (p1) -- (p5);
\draw[->-] (p0) -- (p4);
\draw[->-] (p2) -- (p7);
\draw[->-] (p1) -- (p6);
\draw[->-] (p3) -- (p5);
\draw[->-] (p4) -- (p7);
\draw[->-] (p4) -- (p6);
\draw[->-] (p3) -- (p6);
\draw[->-] (p5) -- (p9);
\draw[->-] (p1) -- (p11);
\draw[->-] (p9) -- (p11);
\draw[->-] (p10) -- (p11);
\draw[->-] (p8) -- (p9);
\draw[->-] (p6) -- (p10);
\draw[->-] (p3) -- (p7);
\draw[->-] (p3) -- (p8);
\draw[->-] (p8) -- (p10);

\draw[fill=white] (p0) circle (2mm) node {0};
\draw[fill=white] (p1) circle (2mm) node {1};
\draw[fill=white] (p2) circle (2mm) node {2};
\draw[fill=white] (p3) circle (2mm) node {3};
\draw[fill=white] (p4) circle (2mm) node {4};
\draw[fill=white] (p5) circle (2mm) node {5};
\draw[fill=white] (p6) circle (2mm) node {6};
\draw[fill=white] (p7) circle (2mm) node {7};
\draw[fill=white] (p8) circle (2mm) node {8};
\draw[fill=white] (p9) circle (2mm) node {9};
\draw[fill=white] (p10) circle (2mm) node {\small10};
\draw[fill=white] (p11) circle (2mm) node {\small11};

\node[align=center] at ($(p0)!0.5!(p1)!0.3333!(p2)+(2mm,1mm)$) {cell 0\\[-3pt]\footnotesize{\pyth{[0, 2, 1, 5,}}\\[-3pt]\footnotesize{\pyth{  4, 7, 6, 3]}}};
\node[align=center] at ($(p1)!0.5!(p11)!0.3333!(p5)+(1mm,1mm)$) {cell 1\\[-3pt]\footnotesize{\pyth{[5, 9, 3, 8,}}\\[-3pt]\footnotesize{\pyth{    1, 11, 6, 10]}}};
\end{3dtikzpicture}
  \subcap{The mesh of hexahedra used in the third example, with
  low-to-high edge orientations.}
  \label{fig:p4_mesh_def}
  \end{subfigure}

  \begin{subfigure}{\textwidth}
  \centering
  \input{img/p4_mesh_numbered_cells}
  \subcap{The local numbering and orientations of the entities of the
  three cells in the third example. The edges labelled with stars are
  those whose local orientation does not match their global orientation
  shown in \cref{fig:p4_mesh_def}.}
  \label{fig:p4_mesh_numbered_cells}
  \end{subfigure}
  \phantomcaption
  \label{fig:p4__mesh}
\end{figure}

\begin{table}
\begin{tabular}{lcccccccccccc}
&\DOFedge{0}&\DOFedge{1}&\DOFedge{2}&\DOFedge{3}&\DOFedge{4}&\DOFedge{5}&\DOFedge{6}&\DOFedge{7}&\DOFedge{8}&\DOFedge{9}&\DOFedge{\tiny10}&\DOFedge{\tiny11}\\[2pt]
cell 0&false&false&false&false&false&false&false&true&false&false&true&true\\[1pt]
cell 1&false&true&true&true&false&false&false&false&false&false&true&false
\end{tabular}

\vspace{2mm}

\begin{tabular}{lcccccc}
  &\DOFface{0}&\DOFface{1}&\DOFface{2}&\DOFface{3}&\DOFface{4}&\DOFface{5}\\[2pt]
  cell 0&true, 0&false, 0&false, 0&false, 0&false, 0&true, 2\\[1pt]
  cell 1&true, 2&true, 2&false, 1&true, 1&true, 0&true, 0
  \end{tabular}
  \vspace{3mm}
  \caption{Do the edges (top) and faces (bottom) of each cell of the mesh
  in \cref{fig:p4__mesh} need reflecting,
  and how many times do the faces need rotating?}
  \label{table:p4_reflect_data}
\end{table}

\subsection{Extension to arbitrary cells and mixed meshes}

\begin{figure}
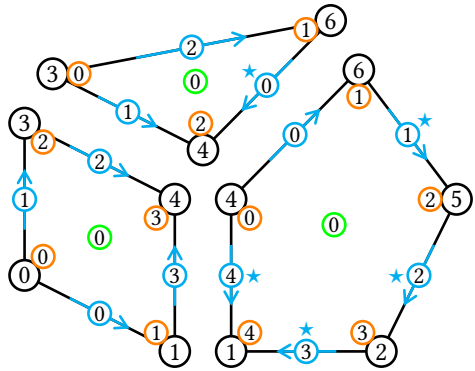

  \newcommand{\nodeDOF}[3]{
\draw[color1, fill=white] ($(#1)+({3.5mm*cos(#2)},{3.5mm*sin(#2)})$) circle (1.5mm) node[black] {\small#3};
 }

\newcommand{\edgeDOF}[4][]{
\draw[color2,->] ($(#2)!0.3!(#3)$) -- ($(#2)!0.7!(#3)$);
\draw[color2, fill=white] ($(#2)!0.5!(#3)$) circle (1.5mm) node[black] {\small#4};
\ifthenelse{\equal{#1}{}}{}{
\node at ($(#2)!0.5!(#3) + ({3mm*cos(#1)},{3mm*sin(#1)})$) {\color{color2}$\star$};
}
}

\begin{2dtikzpicture}
\begin{scope}
\coordinate (p0) at (0,1);
\coordinate (p1) at (2,0);
\coordinate (p3) at (0,3);
\coordinate (p4) at (2,2);

\draw (p0) -- (p1) -- (p4) -- (p3) -- cycle;

\draw[fill=white] (p0) circle (2mm) node[black] {0};
\draw[fill=white] (p1) circle (2mm) node[black] {1};
\draw[fill=white] (p3) circle (2mm) node[black] {3};
\draw[fill=white] (p4) circle (2mm) node[black] {4};

\draw[color3, fill=white] ($(p0)!0.5!(p4)$) circle (1.5mm) node[black] {\small0};
\nodeDOF{p0}{45}{0}
\nodeDOF{p1}{135}{1}
\nodeDOF{p3}{-45}{2}
\nodeDOF{p4}{-135}{3}
\edgeDOF{p0}{p1}{0}
\edgeDOF{p0}{p3}{1}
\edgeDOF{p3}{p4}{2}
\edgeDOF{p1}{p4}{3}
\end{scope}
\begin{scope}[shift={(0.75,0)}]
\coordinate (p1) at (2,0);
\coordinate (p2) at (4,0);
\coordinate (p4) at (2,2);
\coordinate (p5) at (5,2);
\coordinate (p6) at (3.7,3.7);

\draw (p4) -- (p1) -- (p2) -- (p5) -- (p6) -- cycle;

\draw[fill=white] (p4) circle (2mm) node[black] {4};
\draw[fill=white] (p1) circle (2mm) node[black] {1};
\draw[fill=white] (p5) circle (2mm) node[black] {5};
\draw[fill=white] (p6) circle (2mm) node[black] {6};
\draw[fill=white] (p2) circle (2mm) node[black] {2};

\draw[color3, fill=white] ($(p1)!0.55!(p2)!0.45!(p6)$) circle (1.5mm) node[black] {\small0};

\nodeDOF{p4}{-45}{0}
\nodeDOF{p6}{-90}{1}
\nodeDOF{p5}{180}{2}
\nodeDOF{p2}{135}{3}
\nodeDOF{p1}{45}{4}
\edgeDOF{p4}{p6}{0}
\edgeDOF[40]{p6}{p5}{1}
\edgeDOF[180]{p5}{p2}{2}
\edgeDOF[90]{p2}{p1}{3}
\edgeDOF[0]{p4}{p1}{4}
\end{scope}

\begin{scope}[shift={(.375,0.6495)}]
\coordinate (p3) at (0,3);
\coordinate (p4) at (2,2);
\coordinate (p6) at (3.7,3.7);

\draw (p4) -- (p3) -- (p6) -- cycle;

\draw[fill=white] (p3) circle (2mm) node[black] {3};
\draw[fill=white] (p4) circle (2mm) node[black] {4};
\draw[fill=white] (p6) circle (2mm) node[black] {6};

\draw[color3, fill=white] ($(p3)!0.5!(p4)!0.33333!(p6)$) circle (1.5mm) node[black] {\small0};

\nodeDOF{p3}{0}{0}
\nodeDOF{p6}{200}{1}
\nodeDOF{p4}{95}{2}

\edgeDOF[140]{p6}{p4}{0}
\edgeDOF{p3}{p4}{1}
\edgeDOF{p3}{p6}{2}

\end{scope}

\end{2dtikzpicture}
  \caption{Three cells in a mixed mesh. The edges and faces labelled
  with stars are those whose local orientation does not match their
  global orientation.}
  \label{fig:mixed_mesh_numbered_cells}
\end{figure}

The proposed method considers sub-entities of each cell, and the same
transformation is applied to a sub-entity of a given type whatever the type
of the cell. It is straightforward, therefore, to extend this method to
arbitrary polygons and polyhedra and meshes with a mixture of cell
types. \Cref{fig:mixed_mesh_numbered_cells} shows three example cells in
a mixed mesh. As in the previous examples, we could apply the base
transformations for the edges labelled with a star.

On a mesh of hexagonal prisms, for example, each cell would have
eighteen edges and eight faces: one base transformation would be needed
for each edge, and two base transformations would be needed for each
face. As in the previous examples, these base transformations can then
be combined to find the transformation for any cell.

\section{Implementation}
\label{sec:implementation}

The DOF transformation method described in this paper is implemented in
the FEniCSx open source finite element libraries \cite{fenicsbook,
alnaes:2015}. In Basix---the FEniCSx library for element definitions,
basis function computation and interpolation---we compute the base
transformations for each element and provide these alongside the other
information about the finite element defined on the reference cell.
FFCx---the FEniCSx form compiler---generates kernels (in C) from the
Unified Form Language (UFL) \cite{ufl} to evaluate finite element forms
on a cell. Evaluation is on a common reference cell and the kernels are
executed without any reference to cell orientations. DOLFINx is the main
user interface of FEniCSx and executes cell kernels on a mesh. The
transformation of cell basis functions, discussed in
\cref{sec:scalar-spaces,sec:div-spaces,sec:curl-spaces,sec:composing-transformations},
is managed by DOLFINx using Basix: after calling the kernel for a cell,
DOLFINx calls Basix to apply, in place, the required transformations to
the kernel output. DOLFINx passes the required cell orientation data to
Basix. The data can be packed into a single integer, with one bit per
entity used to signify whether or not the entity is reflected, and 2
bits per face (for triangles and quadrilaterals; polygons with more
edges require more bits) to signify how many times the face has been
rotated. The information for each entity can then be retrieved using bit
operations.\label{sec_with_bit_packing} For elements where the base
transformations are simply permutations of the DOFs (for example
Lagrange elements), the effects of the permutations is incorporated into
the DOF map construction and no transformations are required during the
assembly process.

The transformation matrices $M$ for each cell are never explicitly
computed. Instead the necessary composition of base transformations is
directly applied to the evaluations of the reference basis functions. As
each base transformation only affects a small number of basis functions,
we can restrict our transformation to only affect the evaluations of the
relevant functions rather than repeatedly multiplying by $n$ by $n$
matrices.

\subsection{Computing the base transformations in Basix}
\label{serendipity_example}

For Lagrange spaces, it is reasonably straightforward to implement
functions that generate the necessary permutations of the DOF numbering
for an arbitrary degree element. For spaces defined using integral
moments, the effect of permuting a cell sub-entity is dependent upon the
space against which the integral moments are taken. As an example, we
consider a degree 5 serendipity element on a hexahedron
\cite{Arnold_2011}. The functionals \(l_0\), ..., \(l_{73}\) that define
this space are defined by:
\begin{itemize}
\item \(l_0\) to \(l_7\) are point evaluations at the vertices of the
hexahedron.

\item \(l_8\) to \(l_{55}\) are integral moments against a degree 3
Lagrange space on the edges of the hexahedron. There are 4 of these
functionals for each edge.

\item \(l_{56}\) to \(l_{73}\) are integral moments against a degree 1
dPc space on the faces of the hexahedron. A dPc space on a quadrilateral
contains the basis functions of a Lagrange space on a
triangle~\cite{femtable}: the basis functions of a degree 1 dPc space
are:
\begin{align*}
\psi_0(t, s) &= 1-t-s,\\
\psi_1(t, s) &= t,\\
\psi_2(t, s) &= s.
\end{align*}
There are 3 of these functionals for each face.
\end{itemize}
To illustrate how the base transformations can be computed, we consider
the DOFs of this element on the 0th face of the cell ($l_{56}$ to
$l_{58}$). To compute the base transformation corresponding to a
rotation of this face, we consider the transformation $T:(t,s)\mapsto
(s, 1-t)$. This represents a clockwise rotation of the face, and can be
used to undo anticlockwise rotations. If we apply the transformation $T$
to the basis functions of the dPc space on the face, we obtain
\begin{align*}
\psi_0(t, s) &\mapsto \tilde{\psi}_0(t, s) = t-s,\\
\psi_1(t, s) &\mapsto \tilde{\psi}_1(t, s) = s,\\
\psi_2(t, s) &\mapsto \tilde{\psi}_2(t, s) = 1-t.
\end{align*}
These can be represented in terms of $\psi_0$ to $\psi_2$ as
\begin{align}
\begin{pmatrix}
\tilde{\psi}_0\\
\tilde{\psi}_1\\
\tilde{\psi}_2
\end{pmatrix}
&=
\begin{pmatrix}
0&1&-1\\
0&0&1\\
1&0&1
\end{pmatrix}
\begin{pmatrix}
\psi_0\\
\psi_1\\
\psi_2
\end{pmatrix}.\label{the_interpolation}
\end{align}
The matrix here contains the coefficients of an interpolation of the
functions $\tilde{\psi}_0$ to $\tilde{\psi}_2$ into the dPc space on the
quadrilateral face. The same matrix can be use to write the functionals
$\tilde{l}_{56}$ to $\tilde{l}_{58}$ (representing integral moments with
the functions $\tilde{\psi}_0$ to $\tilde{\psi}_2$) in terms of $l_{56}$
to $l_{58}$:
\begin{align*}
\begin{pmatrix}
\tilde{l}_{56}\\
\tilde{l}_{57}\\
\tilde{l}_{58}\\
\end{pmatrix}
&=
\begin{pmatrix}
0&1&-1\\
0&0&1\\
1&0&1
\end{pmatrix}
\begin{pmatrix}
l_{56}\\
l_{57}\\
l_{58}
\end{pmatrix}.
\end{align*}

The basis functions $\phi_{56}$ to $\phi_{58}$ of the serendipity space,
and the basis functions $\tilde{\phi}_{56}$ to $\tilde{\phi}_{58}$
corresponding to the functionals $\tilde{l}_{56}$ to $\tilde{l}_{58}$
satisfy
\begin{align*}
\begin{pmatrix}
\tilde{l}_{56}\\
\tilde{l}_{57}\\
\tilde{l}_{58}\\
\end{pmatrix}
\begin{pmatrix}
\tilde{\phi}_{56}\\
\tilde{\phi}_{57}\\
\tilde{\phi}_{58}\\
\end{pmatrix}^t
&=
\begin{pmatrix}
1&0&0\\
0&1&0\\
0&0&1
\end{pmatrix},
&
\begin{pmatrix}
{l}_{56}\\
{l}_{57}\\
{l}_{58}\\
\end{pmatrix}
\begin{pmatrix}
{\phi}_{56}\\
{\phi}_{57}\\
{\phi}_{58}\\
\end{pmatrix}^t
&=
\begin{pmatrix}
1&0&0\\
0&1&0\\
0&0&1
\end{pmatrix}.
\end{align*}
The base transformation that we want to calculate is the matrix $A$
such that
\begin{align*}
\begin{pmatrix}
\phi_{56}\\
\phi_{57}\\
\phi_{58}
\end{pmatrix}
&=
A
\begin{pmatrix}
\tilde{\phi}_{56}\\
\tilde{\phi}_{57}\\
\tilde{\phi}_{58}\\
\end{pmatrix}.
\end{align*}
Combining the above, we see that
\begin{align*}
\begin{pmatrix}
1&0&0\\
0&1&0\\
0&0&1
\end{pmatrix}
&=
\begin{pmatrix}
\tilde{l}_{56}\\
\tilde{l}_{57}\\
\tilde{l}_{58}\\
\end{pmatrix}
\begin{pmatrix}
\tilde{\phi}_{56}\\
\tilde{\phi}_{57}\\
\tilde{\phi}_{58}\\
\end{pmatrix}^t
\\&=
\begin{pmatrix}
0&1&-1\\
0&0&1\\
1&0&1
\end{pmatrix}
\begin{pmatrix}
l_{56}\\
l_{57}\\
l_{58}
\end{pmatrix}
\begin{pmatrix}
\phi_{56}\\
\phi_{57}\\
\phi_{58}\\
\end{pmatrix}^tA^{-t}
\\&=
\begin{pmatrix}
0&1&-1\\
0&0&1\\
1&0&1
\end{pmatrix}A^{-t},
\end{align*}
and therefore the base transformation representing a rotation of the 0th
face is given by
\begin{align*}
A &=
\begin{pmatrix}
0&1&-1\\
0&0&1\\
1&0&1
\end{pmatrix}^{t}
=\begin{pmatrix}
0&0&1\\
1&0&0\\
-1&1&1
\end{pmatrix}.
\end{align*}

We can therefore compute the base transformations for DOFs defined by
integral moments by computing the interpolation matrix into the space we
are taking integral moments against, as shown in
\cref{the_interpolation}. In Basix, the information needed to calculate
these interpolation coefficients is created when an element is
initialised, so this information can be reused to compute the base
transformations.

The base transformations for the reflection of a face, and for integral
moments on other mesh sub-entities can be obtained by following a
similar method.

\section{Conclusions}
\label{sec:conclusion}

We have presented a method for constructing Ciarlet finite element
spaces of arbitrary degree for any cell type, and for a range of
degrees-of-freedom types, and hence finite element types. The method
does not depend on re-ordering the mesh, which removes the major
limitation of other approaches. The method can be applied to meshes of
triangles, tetrahedra, quadrilaterals and hexahedra, and  it also
extends to general polygonal and polyhedral cells, and mixed meshes.

For some elements, this approach involves applying a permutation to the
DOF numbering on each cell. For other elements, the issues cannot be
resolved using permutations so linear transformations must be applied to
the local element tensor for each cell.

The construction of the entity permutation data needed by this method
only needs knowledge of the global numbering of the vertices in the
mesh. This is an advantage over ordering the mesh, as the input ordering
is preserved, no parallel communication is required, and the method is
valid for all meshes of any cell type. The approach is particularly
appealing in high performance computing applications as it allows the
use of high degree elements on meshes of hexahedral cells.

\begin{acks}
Support from EPSRC (EP/S005072/1) and Rolls-Royce plc as part of the
Strategic Partnership in Computational Science for Advanced Simulation
and Modelling of Engineering Systems (ASiMoV) is gratefully
acknowledged.
\end{acks}
\bibliographystyle{ACM-Reference-Format}
\bibliography{refs}
\end{document}